\documentclass[conference,a4paper,10pt]{IEEEtran}

\usepackage{lmodern}
\usepackage[english]{babel}
\usepackage[utf8]{inputenc} 
\usepackage{amsmath}
\usepackage{amsthm}
\usepackage{hyperref}
\usepackage{enumitem}
\usepackage{xparse}%
\usepackage{xkeyval}%
\usepackage{subfloat}
\usepackage{algorithm}
\usepackage{setspace} 
\usepackage{algorithmic}  
\usepackage[algo2e,noend]{algorithm2e} 
\usepackage{amsfonts}
\usepackage{graphicx}
\usepackage{subfig}
\usepackage{caption}
\usepackage{xcolor}
\usepackage{breqn}
\usepackage{caption}
\usepackage{float}
\usepackage{bm}
\usepackage{natbib}
\usepackage{amssymb}
\usepackage{xspace}
\usepackage[tableposition=top]{caption}
\usepackage{booktabs}
\usepackage{etoolbox}
\usepackage{textcomp}

\newtheorem{modl}{Model}

\newenvironment{model}{\begin{samepage}\begin{modl}}{\end{modl}\end{samepage}}

\newcommand{\qmax}{\overline{q}_i}
\newcommand{\qmin}{\underline{q}_i}

\newcommand{\z}{z_{ij}}

\newcommand{\tm}{\munderbar{t}_{ij}}

\makeatletter
\def\munderbar#1{\underline{\sbox\tw@{$#1$}\dp\tw@\z@\box\tw@}}
\makeatother

\makeatletter
\newcommand{\twocolumncaption}{\@dblarg\@twocolumncaption}
\def\@twocolumncaption[#1]#2{%
	\renewcommand{\@makecaption}[2]{%
		\par\vskip\abovecaptionskip\begingroup\small\rmfamily
		\splittopskip=0pt
		\setbox\@tempboxa=\vbox{
			\@arrayparboxrestore \let \\\@normalcr
			\hsize=.5\hsize \advance\hsize-1em
			\let\\\heading@cr
			\noindent ##1\ ##2\par
		}%
		\vbadness=10000
		\setbox\z@=\vsplit\@tempboxa to .55\ht\@tempboxa
		\setbox\z@=\vtop{\hrule height 0pt \unvbox\z@}
		\setbox\tw@=\vtop{\hrule height 0pt \unvbox\@tempboxa}
		\noindent\box\z@\hfill\box\tw@\par
		\endgroup\vskip \belowcaptionskip
	}%
	\setlength{\abovecaptionskip}{4ex}%
	\caption[#1]{#2}%
}
\makeatother

\begin{document}

\pagestyle{empty} 
\title{A two-stage algorithm for aircraft conflict resolution with trajectory recovery}

\author{
	\IEEEauthorblockN{Fernando H. C. Dias, Stephanie Rahme, David Rey}
	\IEEEauthorblockA{School of  Civil and Environmental Engineering, UNSW Sydney\\
		Sydney, 2052, NSW, Australia \\
		\{f.cunhadias,s.rahme.d.rey\}@unsw.edu.au}
}

\maketitle

\begin{abstract}
\textbf{As air traffic volume is continuously increasing, it has become a priority to improve traffic control algorithms to handle future air travel demand and improve airspace capacity. We address the conflict resolution problem in air traffic control using a novel approach for aircraft collision avoidance with trajectory recovery. We present a two-stage algorithm that first solves all initial conflicts by adjusting aircraft headings and speeds, before identifying the optimal time for aircraft to recover towards their target destination. The collision avoidance stage extends an existing mixed-integer programming formulation to heading control. For the trajectory recovery stage, we introduce a novel exact mixed-integer programming formulation as well as a greedy heuristic algorithm. The proposed two-stage approach guarantees that all trajectories during both the collision avoidance and recovery stages are conflict-free. Numerical results on benchmark problems show that the proposed heuristic for trajectory recovery is competitive while also emphasizing the difficulty of this optimization problem. The proposed approach can be used as a decision-support tool for introducing automation in air traffic control.}
\end{abstract}
\vspace{0.5cm}

\section{Introduction}
Air traffic control (ATC) requires decisions to be quick, effective and free of error. Given the increasing demand observed in the last decades and limited airspace environment, there is an intensifying demand for this service, which is reflected in the increasing amount of control necessary to guarantee safety. Denser and congested traffic configurations may impair flight safety. Yet, state-of-the-art methods for aircraft traffic control are reaching their limits and new approaches, including more automation, have recently received a significant attention in the field \citep{durand2009ant,vela2010near}. Introducing automation within ATC systems has the potential to reduce controller workload and improve airspace capacity \citep{rey2015subliminal}. Conflict detection and resolution (CDR) is vital part of air traffic controllers' workload model. This calls for advanced conflict algorithms capable of acting as efficient decision-support tools. 

The aircraft conflict avoidance and resolution problem can be expressed in the form of an optimization problem, which has the objective to find conflict-free trajectories for all aircraft in a delimited airspace. Many strategies have been proposed to address this problem based on the type of maneuvers (applied separately or in combination) that can be issued to aircraft: speed, heading and/or altitude control. Recently, conflict resolution using global optimization has received a growing attention due to its ability to provide optimal solutions that take into account all traffic within an airspace region and are able to consider the overall state. One of the first global optimization approaches for air conflict resolution was introduced by \cite{pallottino2002conflict} which proposed two formulations: one focusing on speed control and another focusing on heading control and both minimize overall flight time. Subsequent approaches proposed speed control and altitude level-assignment to minimize fuel consumption by metering aircraft at conflict points \citep{vela2010near}. In \cite{vela2009two}, the authors proposed a two-stage stochastic optimization model accounting for wind uncertainty and using speed control. Multi-objective optimisation formulations attempting to balance flight deviation with the total number of maneuvers (velocity, heading and/or altitude change), building on the work of \cite{pallottino2002conflict} were proposed by \citep{alonso2011collision,alonso2014exact}. Subliminal speed control methods which focus on speed control only for conflict resolution has also proven to be a powerful and with low impact in terms of deviation and fuel consumption, although it may fail to resolve all conflicts \citep{rey2015subliminal,cafieri2017maximizing}. More recently, \cite{rey2017complex} proposed a complex number formulation for speed and heading control without any form of discretization.

Despite their potential effectiveness, most efforts in conflict resolution have focused on ensuring collision avoidance, but have overlooked the costs and mechanisms for modeling aircraft's recovery to their original trajectory. This may be critical when conflict resolution is performed using heading control which may significantly deviate aircraft from their initial trajectory, thus possibly increasing flight operating costs. Trajectory recovery has received very little attention in the literature due to the challenging nature of the problem. Meta-heuristics such as genetic algorithms \citep{durand1997optimal} and ant colony algorithms \citep{durand2009ant} have been proposed to find conflict-free solutions that ensure aircraft to avoidance and recovery safely. \citet{dougui2013light} proposed a model which uses an analogy with light propagation theory to create conflict-free aircraft trajectories with recovery. \citet{peyronne2015solving} proposed a B-splines model which uses way-points of a given trajectory to design conflict-free trajectories with recovery. In \cite{omer2015space}, the authors proposed a formulation providing parallel trajectory recovery while minimizing fuel consumption and delays. In this model, aircraft are assumed to perform a preventive maneuver before intersection and the formulation is focused on separating aircraft on their parallel  trajectories. Heading angles are discretised and the optimisation controls both aircraft heading and recovery time. Recently, \citet{lehouillier2017two} proposed a maneuver-discretized model in which pre-defined sets of maneuvers are available for aircraft and a clique-based formulation is proposed to find the optimal combination of confict-free maneuvers. This review of the literature highlights that despite recent improvements in computational optimization, there remain significant open challenges in the design of scalable and exact global optimization approaches for conflict resolution in air traffic control, especially on how to incorporate recovery in a scalable and effective way. 

In this paper, we present a new two-stage algorithm for aircraft conflict resolution with trajectory recovery. In this approach, the speed and heading of aircraft are first optimized to avoid conflicts while minimizing the deviation from their initial trajectories. Then, in a second stage, aircraft trajectories are modified to recover a target position on aircraft's initial trajectories. We next present the mathematical formulation of the proposed conflict resolution model.

\section{Two-stage Conflict Resolution with Trajectory Recovery}

In this section, we present a two-stage approach for conflict resolution with trajectory recovery. We assume that aircraft current and target positions are known and are conflict-free. This sets the context of the optimization problem of interest: given a set of aircraft with known current and target positions, find least-deviating conflict-free trajectories for all aircraft, such that aircraft may safely reach their target destination. To address this problem, we propose to decompose the trajectory optimization problem in two stages: 1) collision avoidance and 2) trajectory recovery. The first stage focuses on controlling aircraft heading and speed to avoid all conflicts while the second stage focuses on calculating the optimal time for aircraft to start safely recovering towards their target position. For brevity, we focus on the two-dimensional conflict resolution problem and only consider horizontal aircraft maneuvers. The extension to the vertical case can be addressed by incorporating flight level change maneuvers in the collision avoidance stage \citep{dias2019disjunctive} and ensuring safe recovery to aircraft target flight level and position. We leave this extension for future research. 

\subsection{Collision Avoidance}

In this first stage, the goal is to find conflict-free, least-deviating heading angles and speed changes. 

\subsubsection{Separation Conditions}

Consider a set of aircraft $\mathcal{A}$ sharing the same flight level. For each aircraft $i \in \mathcal{A}$, assuming uniform motion laws, its position is: $p_i(t) = [x_i(t) = \widehat{x}_i + q_iv_i\cos (\widehat{\theta}_{i} + \theta_i)t, y_i(t) = \widehat{y}_i + q_iv_i\sin(\widehat{\theta}_{i} + \theta_i)t]^\top$ in which $v_i$ is the speed, $\widehat{x}_i$ and $\widehat{y}_i$ are the initial coordinates of $i$ at the beginning of its trajectory, $\widehat{\theta}_{i}$ is its initial heading angle, $\theta_i$ is its deviation angle and $q_i$ is the speed deviation. 

To avoid trigonometric functions, we discretize the set of heading change maneuvers. Let $\mathcal{H}_i$ be the set of deviation angles for each aircraft $i \in \mathcal{A}$, and let $\delta_{ik}$ be a binary variable which is 1 if aircraft $i$ selects deviation angle $\theta_k \in \mathcal{H}_i$. The relative velocity vector of $i$ and $j$, denoted $v_{ij}$, can be expressed as $v_{ij} = [v_{ij,x},v_{ij,y}]^\top$ where:
\begin{subequations}
\begin{align}
v_{ij,x} =& q_iv_i\cos\Big(\widehat{\theta_i} + \sum\limits_{k \in \mathcal{H}_i}\delta_{ik}\theta_k\Big) \nonumber \\ 
&- q_jv_j\cos\Big(\widehat{\theta_j} + \sum\limits_{k \in \mathcal{H}_i}\delta_{jk}\theta_k\Big) \\
v_{ij,y} =& q_iv_i\sin\Big(\widehat{\theta_i} + \sum\limits_{k \in \mathcal{H}_i}\delta_{ik}\theta_k\Big) \nonumber \\ 
&- q_jv_j\sin\Big(\widehat{\theta_j} + \sum\limits_{k \in \mathcal{H}_i}\delta_{jk}\theta_k\Big)
\end{align}\label{motion}
\end{subequations}
To linearize the bilinear terms of the form $q_i \delta_{ik}$, we first expand the trigonometric functions using traditional identities and introduce an auxiliary variable $\phi_{ik}\equiv q_i \delta_{ik}$ via the following constraints:
\begin{subequations}
\begin{align}
	\qmin\delta_{ik}\leq \phi_{ik}  \quad  \forall i \in \mathcal{A}, k \in \mathcal{H}_i\\
	\phi_{ik} \leq \delta_{ik}\qmax  \quad  \forall i \in \mathcal{A}, k \in \mathcal{H}_i\\
	q_i - (1 - \delta_{ik})\qmax \leq \phi_{ik}  \quad  \forall i \in \mathcal{A}, k \in \mathcal{H}_i\\
	\phi_{ik} \leq q_i - (1 - \delta_{ik})\qmin  \quad  \forall i \in \mathcal{A}, k \in \mathcal{H}_i
\end{align}\label{phi}
\end{subequations}
The relative position of aircraft $i$ and $j$ at time $t$ can be represented as $p_{ij}(t) = p_{i}(t) - p_{j}(t)$. Let $d=5 NM$ be the horizontal separation norm, two aircraft $i,j \in \mathcal{A}$ are horizontally separated if and only if: $||p_{ij}(t)|| \ge d, \quad \forall t \ge 0$. To derive a time-independent separation condition, we use the method described in \cite{cafieri2017maximizing, cafieri2017mixed, rey2017complex}. Let $f_{ij}(t) = ||v_{ij}(t)||^2t^2 + 2p_{ij}(t)\cdot v_{ij}(t) + ||2p_{ij}||^2 - d^2$ and let $\tm =  \frac{- p_{ij} \cdot v_{ij}}{||v_{ij}||^2}$ be the time at which $f_{ij}(t)$ is minimized, where $v_{ij}$ is the relative velocity vector. Assuming aircraft are initially separated, the separation condition can be represented as $f_{ij}(t) \ge 0$ if $\tm \ge 0$; otherwise aircraft are diverging thus separated. Evaluating $f_{ij}(\tm)$ gives the expression $g_{ij}$ which is time-invariant:
\begin{equation}\label{g_function}
 	g_{ij} = (v_{ij}^x)^2(y_{ij}^2 - d^2) + (v_{ij}^y)^2(x_{ij}^2 - d^2) - v_{ij}^x v_{ij}^y (2 x_{ij} y_{ij}) 
\end{equation}
To distinguish between past and future events, we track aircraft convergence using the sign of $\tm$, which corresponds to the time of minimum separation. Aircraft pairwise separation condition can thus be summarized as:
\begin{equation}\label{sepconditions}
g_{ij} \geq 0 \vee \tm \leq 0
\end{equation}
Following the approach of \citet{rey2017complex}, by isolating variables $v_{ij,x}$ and $v_{ij,y}$ in Eq. \eqref{g_function} and introducing a binary variable $\z$, we obtain the set of constraints:
\begin{subequations}
	\begin{align}
	v_{ij,x}\hat{x}_{ij} - v_{ij,y}\hat{y}_{ij} \leq 0 &\text{ if } \z=1  \\
	v_{ij,x}\hat{x}_{ij} - v_{ij,y}\hat{y}_{ij} \geq 0 &\text{ if } \z=0  \\
	v_{ij,x} \psi_{ij}^l - v_{ij,y} \varphi_{ij}^l \leq 0 &\text{ if } \z=1 \\
	v_{ij,x} \psi_{ij}^u - v_{ij,y} \varphi_{ij}^u \geq 0 &\text{ if } \z=0 
	\end{align}
	\label{eq:sepdis}
\end{subequations}
The coefficients $\psi_{ij}^l$, $\varphi_{ij}^l$ and $\psi_{ij}^u$, $\varphi_{ij}^u$ are constants that can be pre-processed based on the sign of $\hat{x}_{ij}$ and $\hat{y}_{ij}$. \citet{dias2019disjunctive} have shown that imposing these constraints is equivalent to the separation condition Eq. \eqref{sepconditions}.

\subsubsection{Speed Control, Heading Changes and Objective Function}

For each aircraft $i \in \mathcal{A}$, we assume that the speed rate variable is lower bounded by $\qmin$ and upper bounded by $\qmax$, thus the speed control constraint is:
\begin{equation}
\qmin \leq q_i \leq \qmax \qquad \forall i \in \mathcal{A}
\label{speedbound}
\end{equation}
To model heading angle changes, we assume that each aircraft $i \in \mathcal{A}$ has access to a set of options for heading angles changes $k \in \mathcal{H}_i$. The selection is given by the binary variable $\delta_{ik}$ which is equal to 1 if aircraft $i$ selects and angle $\theta_k$. Heading angle selection is ensured via the constraint:
\begin{equation}
\sum_{k \in \mathcal{H}_i} \delta_{ik} = 1 \qquad \forall i \in \mathcal{A}
\label{angleselection}
\end{equation}
For the objective function, we use a quadratic penalty on speed and heading deviations and the parameter $w \geq 0$ is used to compromise between speed and angle change:
\begin{align}
\text{minimize} &\sum_{i \in \mathcal{A}} \left(w(1 - q_i)^2  + (1-w)\left(\sum_{k \in \mathcal{H}_i} \delta_{ik}\theta_{k}\right)^2\right) \label{quadobj}
\end{align}

\subsubsection{Disjunctive formulation}

The proposed approach for conflict resolution via speed and discrete heading control is summarized in Model \ref{mod:dm}. 
\begin{model}[Collision Avoidance]
	\label{mod:dm}
	\allowdisplaybreaks
	\begin{subequations}
		\begin{align*}
		&\text{\emph{minimize}}\qquad \eqref{quadobj}\\
		&\text{\emph{subject to}}  && \nonumber \\
		& \eqref{motion},\eqref{phi},\eqref{eq:sepdis}, \eqref{speedbound},\eqref{angleselection} \\
		& v_{ij,x},v_{ij,y} \in \mathbb{R} &&\forall (i,j) \in \mathcal{P} \\  
		& \z \in \{0,1\} &&\forall (i,j) \in \mathcal{P} \\
		& \phi_{ik} \in \{0,1\} && i \in \mathcal{A}, k \in \mathcal{H}_i
		\end{align*}
	\end{subequations}
\end{model}

This formulation is a mixed-integer quadratic program (MIQP). The model is built with indicator constraints which can be handled by commercial optimization software. In our implementation, we handle indicator constraints by deriving the convex hull of these constraints using the method proposed by \cite{hijazi2010mixed}. For implementation details, a fully reproducible model can be found at: \small\url{https://github.com/davidrey123/Conflict\_Resolution\_for\_Air\_Traffic_Control}\normalsize. 

The outputs of Model \ref{mod:dm} are vectors of heading angle speed deviation, ${\delta}^\star$ and ${q}^\star$. Let ${\theta_{i}^\star}$ be the optimal heading change for aircraft $i$, we have: ${\theta_{i}^\star} = \sum\limits_{k \in \mathcal{H}_i} {\delta_{ik}^\star}\theta_{k}$. We use the resulting conflict-free trajectories of Model \ref{mod:dm} as input for the recovery stage of the proposed conflict resolution algorithm. 

\subsection{Trajectory Recovery}

\begin{figure}
	\centering
	\includegraphics[height=0.3\linewidth,width=0.7\linewidth]{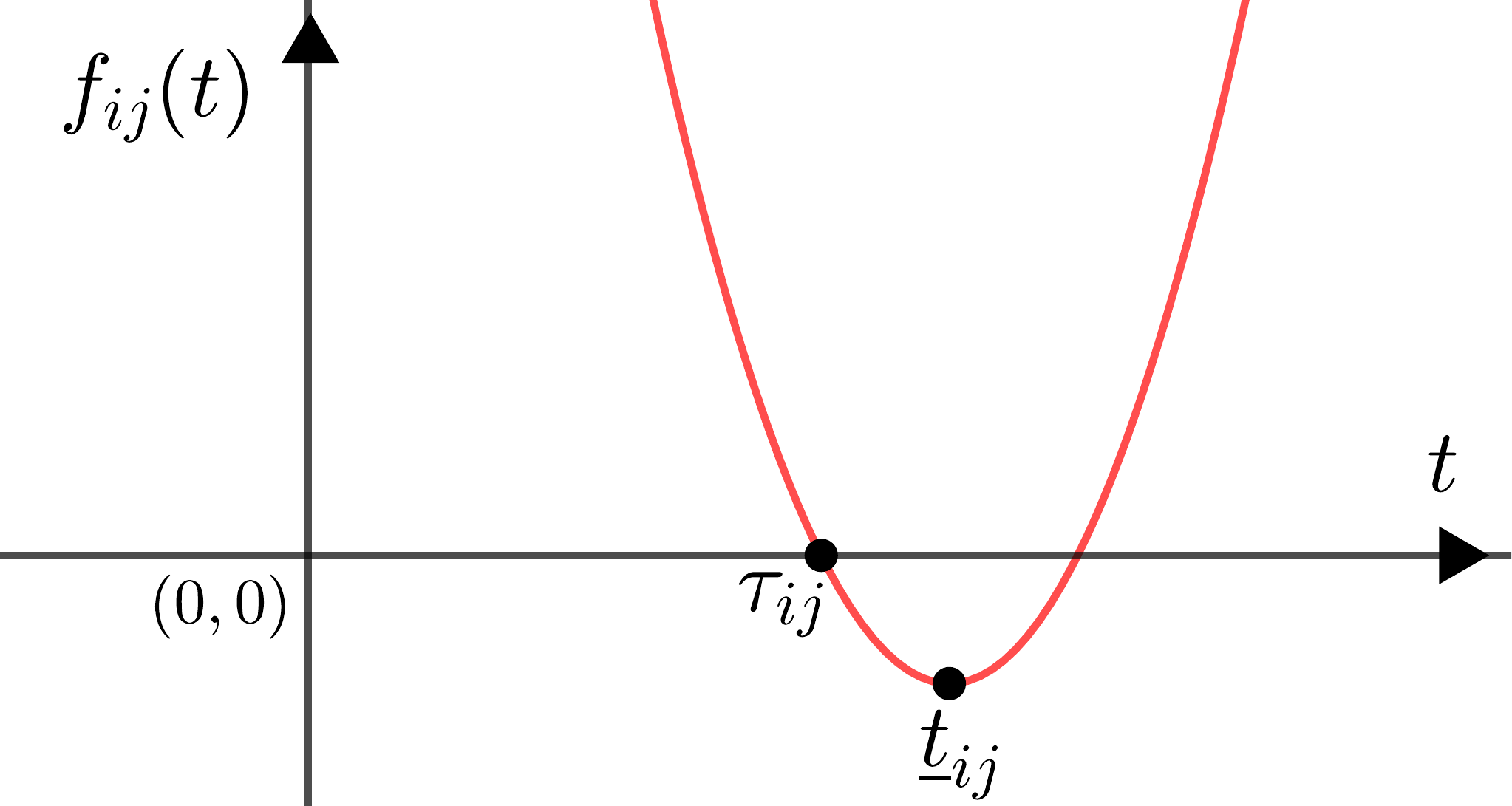}
	\caption{Illustration of $f_{ij}(t)$ for a configuration with $g_{ij} < 0$ and $\tm > 0$. $\tau_{ij}$ represents the start time of the conflict.}
	\label{fig:conflictandturn}
\end{figure}

The second stage aims to identify the optimal time for aircraft to recover towards its target position. To account for the cost of trajectory deviations at the first stage, we denote $a_i^\star$ for $i \in \mathcal{A}$, the deviation corresponding to the optimal solution of Model \ref{mod:dm}:
\begin{equation}\label{eqarea}
a_i^\star = w(1 - q_i^\star)^2+(1-w)(\theta_{i}^\star)^2                               
\end{equation}                                                           Let $t_i$ be the recovery time of aircraft $i \in \mathcal{A}$ and let $\check{x_i},\check{y_i}$ be the coordinates of the target position of  $i$. For trajectory recovery, each aircraft need to perform opposing maneuvers to cancel the deviation applied during avoidance. Similar to the avoidance model, our goal is to guarantee that all pair of aircraft are separated throughout the recovery stage. Since the separation condition Eq. \eqref{sepconditions} is based on linear motion, we need to distinguish the trajectory stage of each aircraft $i \in \mathcal{A}$, i.e. before and after its recovery time $t_i$. We denote $A_i$ the avoidance trajectory of aircraft $i$ and $R_i$ its recovery trajectory. Given a pair $(i,j)$ of aircraft, we need to ensure that aircraft are separated during all pairwise trajectory stages, denoted $A_iA_j$, $A_iR_j$, $R_iA_j$ and $R_iR_j$. Observe that separation for the stage $A_iA_j$ is already ensured by the solution of Model \ref{mod:dm}. If aircraft $i$ and $j$ were to recover at the same time period, then aircraft will transition from $A_iA_j$ to $R_iR_j$ directly. Otherwise, if $i$ (resp. $j$) recovers before $j$ (resp. $i$), then $A_iA_j$ will transition to $R_iA_j$ (resp. $A_iR_j$) before transitioning to $R_iR_j$. 

The distance flown during the collision avoidance stage is $d_{A_i}(t_i) = \sqrt{(\widehat{x}_i - x(t_i))^2 + (\widehat{y}_i - y(t_i))^2}$; similarly, the distance flown during the trajectory recovery stage is $d_{R_i}(t_i) = \sqrt{(x(t_i) - \check{x_i})^2 + (y(t_i) - \check{y_i})^2}$. If aircraft $i$ has changed its speed only, then at $t_i$ aircraft should recover its initial speed without any heading change. Otherwise, if $i$ has made a turning movement, then aircraft should turn in the opposite direction at time $t_i$ with the angle:
\begin{equation}
\theta_{R_i}(t_i) = \arcsin\Big(\frac{d_{A_i}(t_i)\sin(\theta_{A_{i}})}{d_{R_i}(t_i)}\Big)
\end{equation}
To avoid trigonometric functions and obtain a tractable formulation, we discretize time, i.e. $t_i \in \{0,1\epsilon, 2\epsilon, \ldots, |\mathcal{T}|\epsilon\}$ where $\mathcal{T}$ is the set of time periods available for recovery and $\epsilon$ is the length of time periods. Abusing notation, we redefine the separation condition expressed in Eq. \eqref{sepconditions} as: $g_{ij}(m,n) \geq 0$ and $\tm(m,n) \leq 0$ where the pair $(m,n)$ indicates the time period indices of recovery times $t_i$ and $t_j$, respectively. 

Let $\Omega_{X_iX_j}$ be the set of conflict-free pairs of recovery times for aircraft $i,j \in \mathcal{A}$ where $X_i$ represents the state of the trajectory of aircraft $i$, \emph{i.e.} $A_i$ or $R_i$; and $Y_j$ represents the state of the trajectory of aircraft $j$, \emph{ i.e.} $A_j$ or $R_j$. This set can be specified into three different sets corresponding to the three different states during the recovery stage. The set $\Omega_{R_iR_j}$ is defined as:
\begin{align}
\Omega_{R_iR_j} &= \{(m,n)\in \mathcal{T}^2 :  g_{R_iR_j}(mn) \geq 0 \vee \munderbar{t}_{R_iR_j}(m,n)\leq 0\}
\end{align}
For the states $A_iR_j$ and $R_iA_j$ an extra condition is required. Consider the state $A_iR_j$: if the lines of motion corresponding to trajectories $A_i$ and $R_j$ are in conflict but aircraft $i$ turns into recovery prior to the start of this conflict, then no conflict will occur. This illustrated in Figure \ref{fig:conflictandturn} where $g_{A_iR_j} < 0$ and $t_{A_iR_j} > 0$. Let $\tau_{A_iR_j}(t_j)$ be the smallest root of $g_{A_iR_j} = 0$ if $j$ recovers at time $t_j$. If aircraft $i$ recovers prior to $\tau_{A_iR_j}(t_j)$, \textit{i.e.} $t_i \leq \tau_{A_iR_j}(t_j)$, then the conflict will be avoided. Accordingly, we define: 
\begin{subequations}
\begin{align}
\Omega_{A_iR_j} &= \{(m,n) \in \mathcal{T}^2 :  g_{A_iR_j}(n) \geq 0 \nonumber \\	
&\quad \vee \munderbar{t}_{A_iR_j}(n)\leq 0 \quad \vee \quad m \leq \tau_{A_iR_j}(n)\} \\
\Omega_{R_iA_j} &= \{(m,n) \in \mathcal{T}^2 :  g_{R_iA_j}(m) \geq 0 \nonumber \\
&\quad \vee \munderbar{t}_{R_iA_j}(m)\leq 0 \quad \vee \quad n \leq \tau_{R_iA_j}(m)\}
\end{align}
\end{subequations}

We next propose an exact and a heuristic approach to optimize aircraft recovery times.

\subsubsection{Exact-Recovery}

Let $\rho_{im}$ be a binary variable equal to 1 if aircraft $i \in \mathcal{A}$ recovers at time period $m \in \mathcal{T}$ and 0 otherwise. We seek to minimize the total weighted recovery time, i.e. $\sum_{i \in \mathcal{A}} \sum_{m  \in \mathcal{T}} a_i^\star\rho_{im}t_m^2$. To track the states of aircraft pair $(i,j)$ which are activated, we introduce two binary variables $\alpha_{ij}$ and $\beta_{ij}$. Those variables are used to identify whether $t_i < t_j$ ($\alpha_{ij} = 1$) which activates state $R_iA_j$, or if $t_i > t_j$ ($\beta_{ij} = 1$) which activates state $A_iR_j$. Variables $\alpha_{ij}$ and $\beta_{ij}$ are defined via the constraints:
\begin{subequations}\label{cons1}
	\begin{align}
	& \alpha_{ij} \geq \frac{1}{|\mathcal{T}|}\Bigg(\sum\limits_{m \in \mathcal{T}}m\rho_{im} - \sum\limits_{n \in \mathcal{T}}n\rho_{jn}\Bigg) &&\forall (i,j) \in \mathcal{P} \\
	& \beta_{ij} \geq \frac{1}{|\mathcal{T}|} \Bigg(\sum\limits_{n \in \mathcal{T}}n\rho_{jn} - \sum\limits_{m \in \mathcal{T}}m\rho_{im} \Bigg) &&\forall (i,j) \in \mathcal{P} \\
	& \alpha_{ij} + \beta_{ij} \leq 1 &&\forall (i,j) \in \mathcal{P}
	\end{align}
\end{subequations}

We use the following constraints to exclude conflicting trajectories from the solution. Observe that states $A_iR_j$ and $R_iA_j$ are conditional on the recovery times of $t_i$ and $t_j$ and thus the corresponding constraints are only active if $i$ and $j$ do not recover at the time period.
\begin{subequations}\label{cons2}
	\begin{align}
	& \rho_{im} + \rho_{jn} \leq 2 - \beta_{ij} &&\forall (i,j) \in \mathcal{P}, (m,n) \in \Omega_{A_iR_j}\\
	& \rho_{im} + \rho_{jn} \leq 2 - \alpha_{ij} &&\forall (i,j) \in \mathcal{P}, (m,n) \in \Omega_{R_iA_j} \\
	& \rho_{im} + \rho_{jn} \leq 1  &&\forall (i,j) \in \mathcal{P}, (m,n) \in \Omega_{R_iR_j}
	\end{align}
\end{subequations}
Aircraft are assigned a recovery time via the constraint:
\begin{align}
& \sum_{m  \in \mathcal{T}} \rho_{im} = 1 \quad \quad \forall i \in \mathcal{A}
\label{cons3}
\end{align}
The exact trajectory recovery formulation is summarized in Model \ref{mod:onestepdetailed} which is a MILP.
\begin{model}[Exact-Recovery]
	\label{mod:onestepdetailed}
	\allowdisplaybreaks
	\begin{align}
	&\text{\emph{minimize}}\qquad \sum_{i \in \mathcal{A}} \sum_{m  \in \mathcal{T}} a_i^\star\rho_{im}t_m^2\\
	&\text{\emph{subject to}} \nonumber \qquad \\
	& \eqref{cons1}, \eqref{cons2}, \eqref{cons3} \nonumber\\
	& \rho_{im} \in \{0,1\}  \hspace{3.8cm} \forall i \in \mathcal{A}, m \in \mathcal{T} \\
	& \alpha_{ij},\beta_{ij} \in \{0,1\}  \hspace{3.78cm} \forall (i,j) \in \mathcal{P}\
	\end{align}
\end{model}

\subsubsection{Greedy-Recovery}

This heuristic iterates over all time steps and uses a priority list to decide which aircraft can be recovered at each time step. The priority list used is based on $a_i^\star$ values \eqref{eqarea}. The algorithm first sorts aircraft accordingly and iterates over time periods. At each time period the algorithm iterates over the sorted list of aircraft and check if each aircraft can be recovered at the current time. The process is repeated until no aircraft can recover at the current time. The proposed algorithm has a worst-case time complexity of $\mathcal{O}(|T||A|^3)$. The pseudo-code of the proposed greedy algorithm for aircraft trajectory recovery is summarized in \ref{algo:greedy}.

\begin{algorithm}[!h]
\setstretch{1.0}
\textbf{Input}: $\mathcal{A}$, $\bm{a}^\star$ \\
\textbf{Output}: $\bm{t}$ \\
$ \mathcal{R} \gets \{i \in \mathcal{A}: a_i = 0\} $ \\
$ \mathcal{D} \gets \mathcal{R} - \mathcal{A} $ \\
$ \mathcal{D} \gets \text{Sort based on decreasing } a_i$ values\\
\For{$t \in \mathcal{T}$}{
	$ \text{update} \gets \texttt{true} $ \\
	\While{$\emph{\text{update}} = \emph{\texttt{true}} $}{
		$ \text{update} \gets \texttt{false} $ \\
		\For{$i \in \mathcal{D}$}{
			$ sep \gets 0$ \\
			\For{$j \in \mathcal{A}$}{
				\If{$i <j$}{
					\eIf{$j \in \mathcal{R}$}{
						\If{$(t,t_j) \in \Omega_{R_iR_j}$}{
							$ sep \gets sep + 1$ 
						}
					}
					{
						\If{$(t,t_j) \in \Omega_{R_iA_j}$}{
							$ sep \gets sep + 1$ 
						}
					}
				}
				\If{$i > j$}{
					\eIf{$j \in \mathcal{R}$}{
						\If{$(t_j,t) \in \Omega_{R_jR_i}$}{
							$ sep \gets sep + 1$ 
						}
					}
					{
						\If{$(t_j,t) \in \Omega_{A_jR_i}$}{
							$ sep \gets sep + 1$
						}
					}
				}
			}
			\If{$sep = |\mathcal{A}|-1$}{$t_i \gets t$\\
			$\mathcal{R} \gets \mathcal{R} \cup \{i\} $\\
			$\mathcal{D} \gets \mathcal{R} \setminus \{i\}$\\
			$\text{update} \gets \texttt{true} $}
		}
	}
}
\caption{Greedy-Recovery Algorithm}
\label{algo:greedy}
\end{algorithm}

We now formally introduce the proposed the two-stage algorithm: stage 1 solves initial conflicts using Model \ref{mod:dm} by adjusting aircraft headings and speeds. The optimal solution of stage 1 is used as input for stage 2 which finds optimal aircraft recovery times. Stage 2 is solved either exactly via Model \ref{mod:onestepdetailed} or using the heuristic algorithm \ref{algo:greedy}. We next present numerical results. 

\section{Numerical Experiments}
\subsection{Experimental Framework}

\begin{figure}[!h]	
	\centering	
	\mbox{\subfloat[CP with 7 aircraft]{\includegraphics[width=4cm]{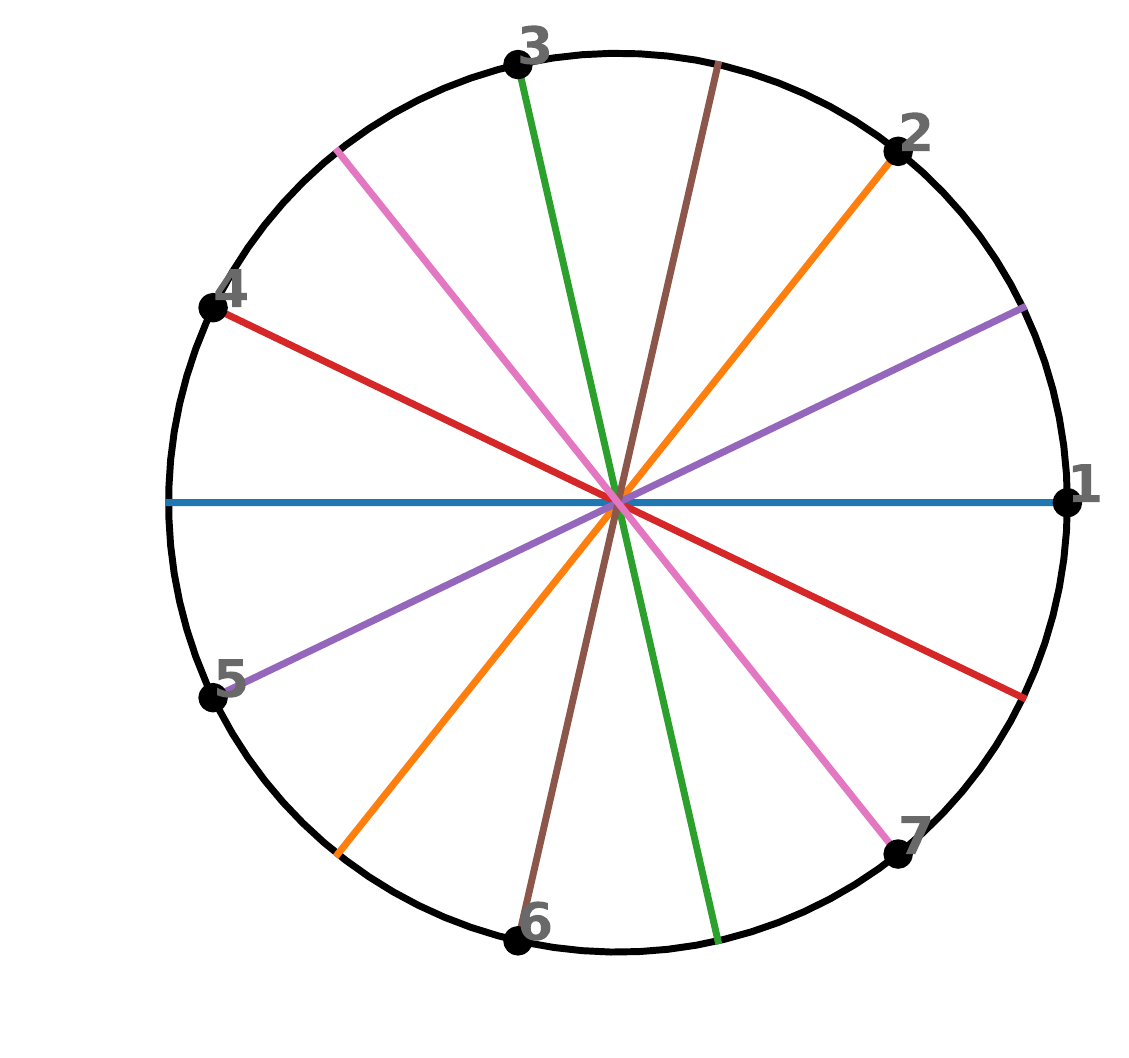}\label{fig:cp}}
		\subfloat[RCP with 10 aircraft]{\includegraphics[width=4cm]{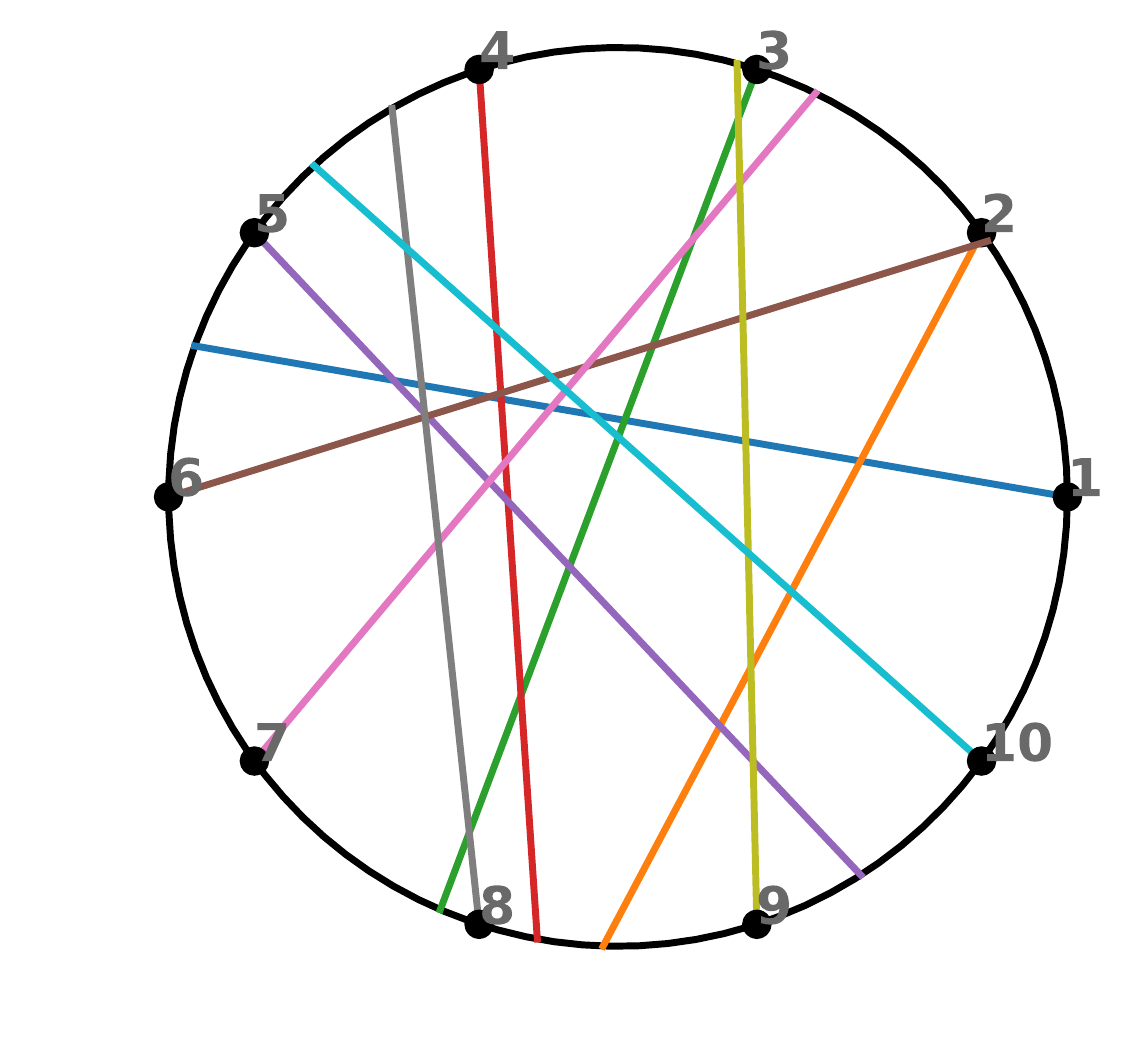} \label{fig:rcp}} } 
	\caption{Example of test instances for the Circle Problem (CP) and Random Circle Problem (RCP)}
	\label{case1}
\end{figure}

We test the performance of the proposed two-stage approach using classical benchmark problems: the Circle Problem (CP) and the Random Circle Problem (RCP). These benchmark problems have been widely used in the community for the horizontal aircraft conflict resolution problem to assess the performance of conflict resolution algorithms \citep{durand2009ant,rey2015equity,cafieri2017maximizing,rey2017complex}. Instances for the CP and RCP are illustrated in Figure \ref{case1}. The CP consists of a set of aircraft uniformly positioned on the circle heading towards its centre. Aircraft speeds are assumed to be identical, hence the problem is highly symmetric (see Fig. \ref{fig:cp}). In contrast, the RCP builds on the same framework, but aircraft initial speeds and headings are randomly deviated within specified ranges to create random instances with less structure (see Fig. \ref{fig:rcp}). 

We use a circle of radius of 200NM for all tests. For CP instances, all aircraft have the same initial speed of 500NM/h. For RCP instances, aircraft initial speeds are randomly chosen in the range 486-594NM/h and their initial headings are deviated from a radial trajectory (i.e. towards the centre of the circle) by adding a randomly chosen an angle between $-\frac{\pi}{6}$ and $+\frac{\pi}{6}$. We report results for CP instances with 4 to 15 aircraft. For RCP instances, we report average performance for three instance size, i.e. with 10, 20 and 30 aircraft. For each instance size, 100 randomly generated RCP instances are generated and solved. All instances are available at: \small\url{https://github.com/davidrey123/Conflict\_Resolution\_for\_Air\_Traffic_Control}\normalsize. 

For all tests, we use a speed regulation range of $\pm10\%$ and allow heading changes in the range $\pm \frac{\pi}{6}$ in steps of $10^{\circ}$, hence a total of 7 headings are available per aircraft (including the initial trajectory heading). We use $w=0.2$ in the objective of Model \ref{mod:dm}. This value was selected such that both heading and speed control terms were of comparable order of magnitude with an emphasis on penalizing heading control. For stage 2, we use a total of $|\mathcal{T}|=15$ time periods, with a step of $\epsilon=2$ minutes. Models \ref{mod:dm} and \ref{mod:onestepdetailed} are solved with CPLEX's Python API and a time limit of 5 minutes.

\subsection{Illustration}
\label{num}

To illustrate the proposed two-stage algorithm, we plot the optimal solution obtained using Model \ref{mod:dm} and Model \ref{mod:onestepdetailed} (ER) for CP instances with 5, 10 and 15 aircraft (CP-5, CP-10 and CP-15). For RCP instances, we show three instances of each instance size tested, i.e. with 10, 20 and 30 aircraft. In the figures, dashed gray lines represent aircraft initial trajectories, red lines represent the avoidance trajectory of stage 1, and blue lines represent recovery trajectories of stage 2. For CP-5, all conflicts are solved using speed control only. Instead, for CP-10 and CP-15, some aircraft make a turn before recovering to their destination. The solutions of RCP instances highlight that increasing the number of aircraft tends to increase the duration of the collision avoidance trajectory (red line). 

\begin{figure*}[!h]
	\centering
	\subfloat[CP-5]{\includegraphics[width=0.2\textwidth]{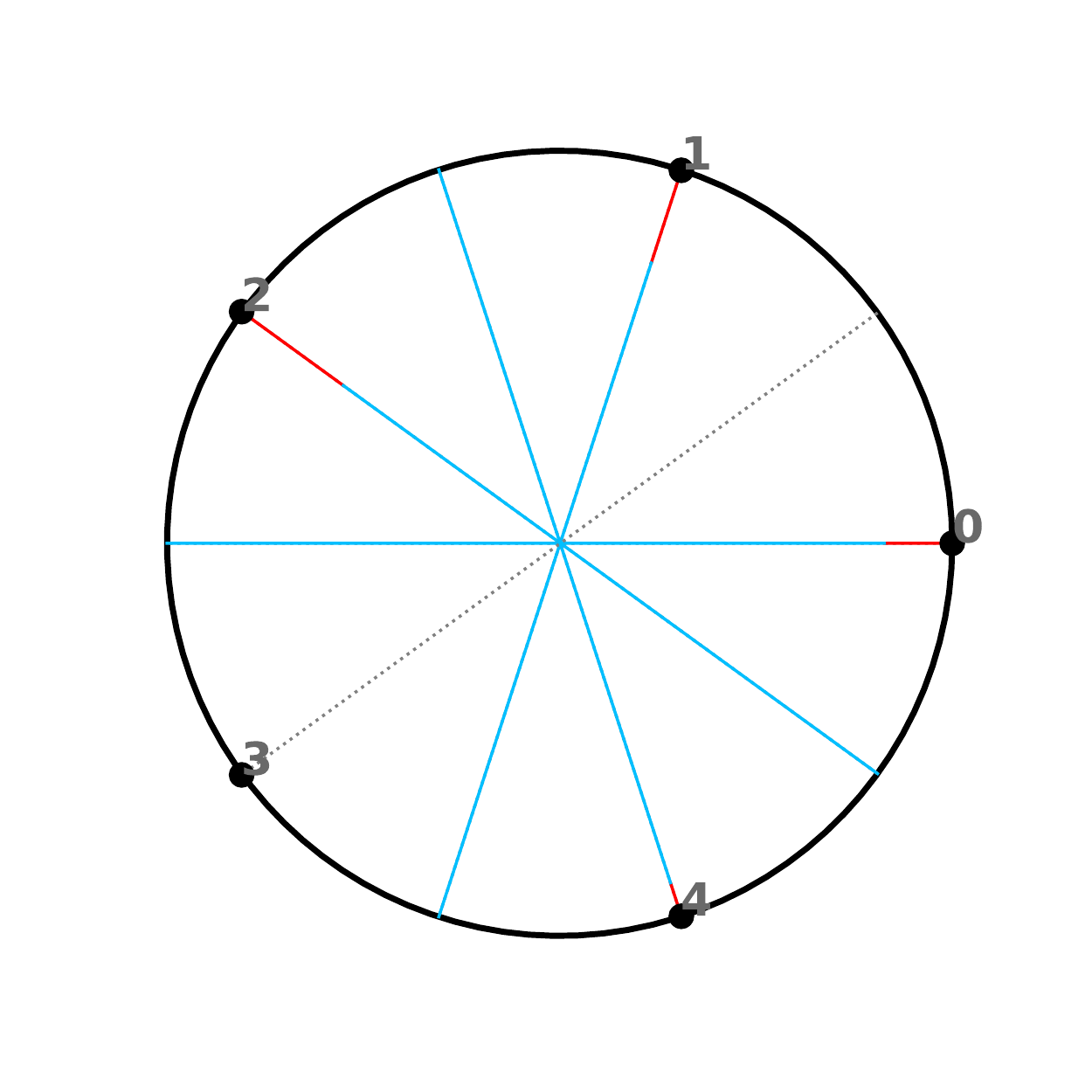}\label{cp1}}\ 
	\subfloat[RCP-10-1]{\includegraphics[width=0.2\textwidth]{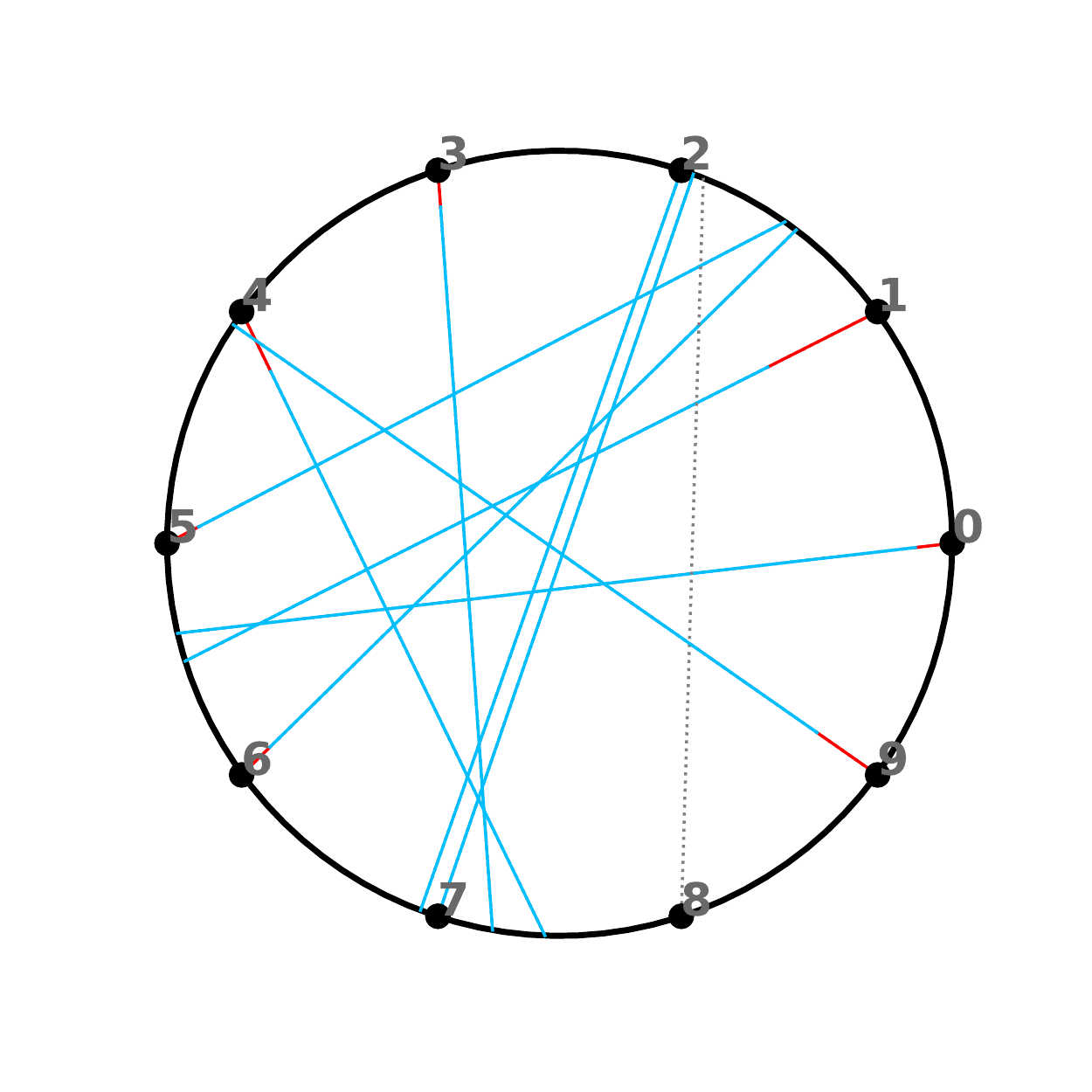}\label{rcp1}}\
	\subfloat[RCP-20-1]{\includegraphics[width=0.2\textwidth]{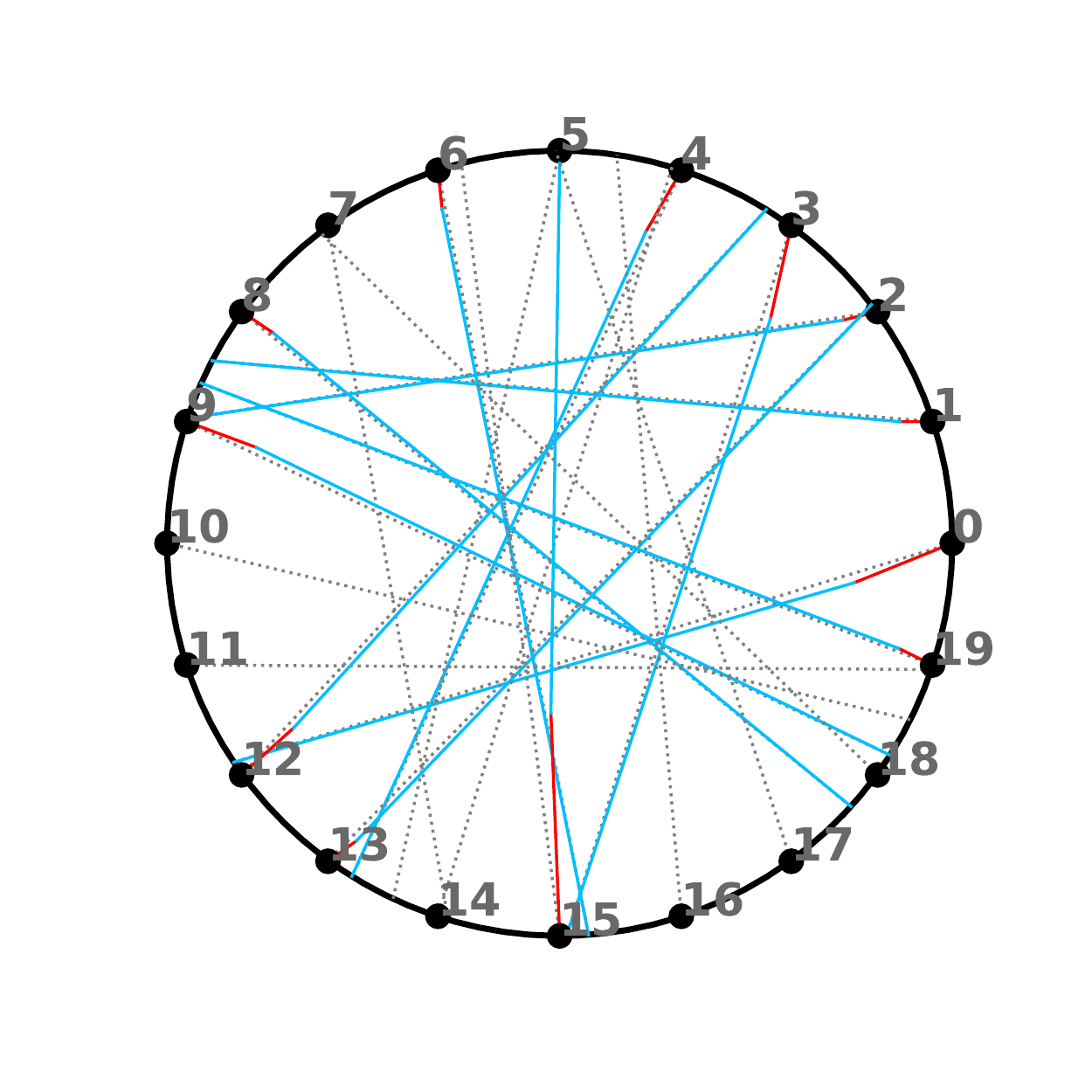}\label{rcp4}}\ 
	\subfloat[RCP-30-1]{\includegraphics[width=0.2\textwidth]{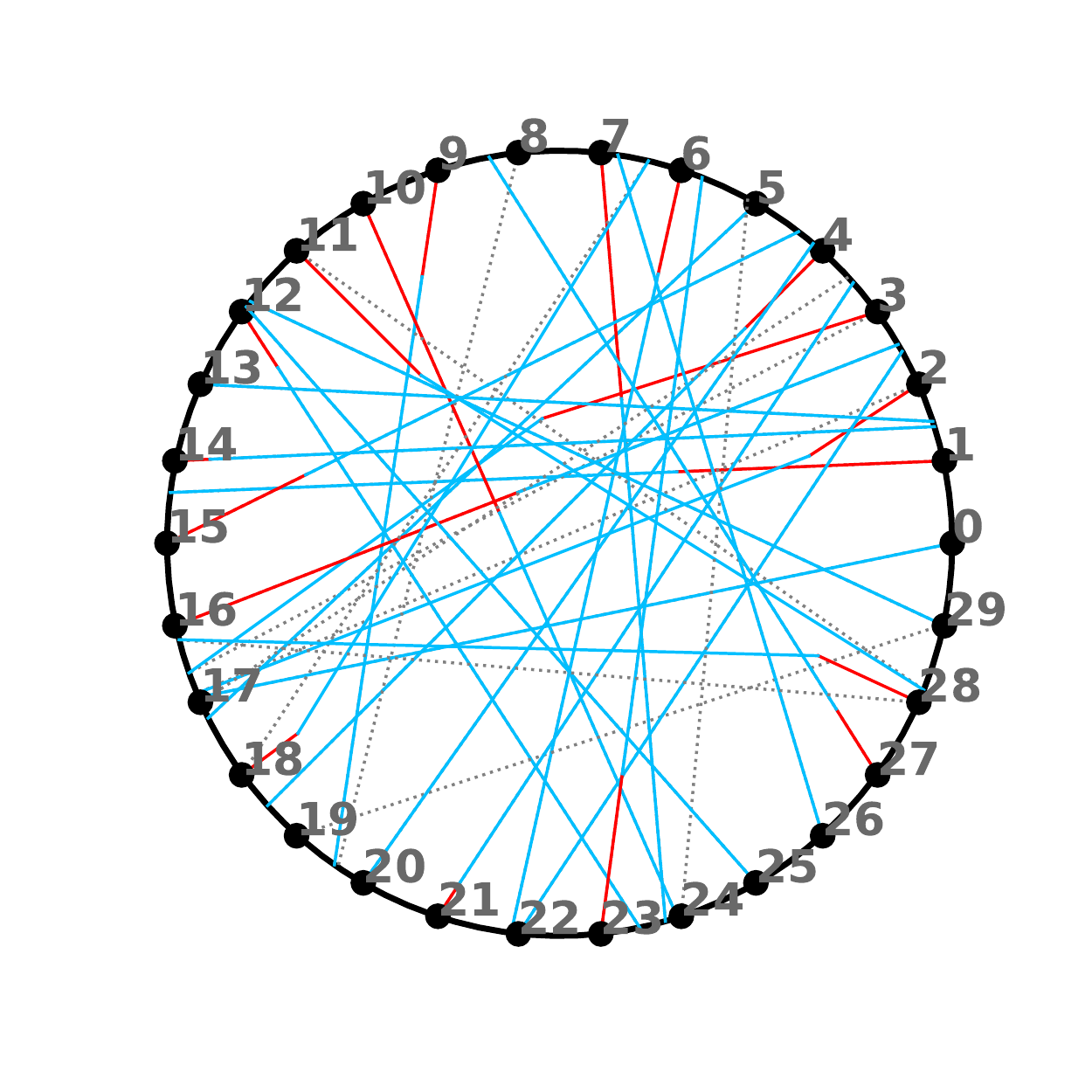}\label{rcp2}}\ 
	\subfloat[CP-10]{\includegraphics[width=0.2\textwidth]{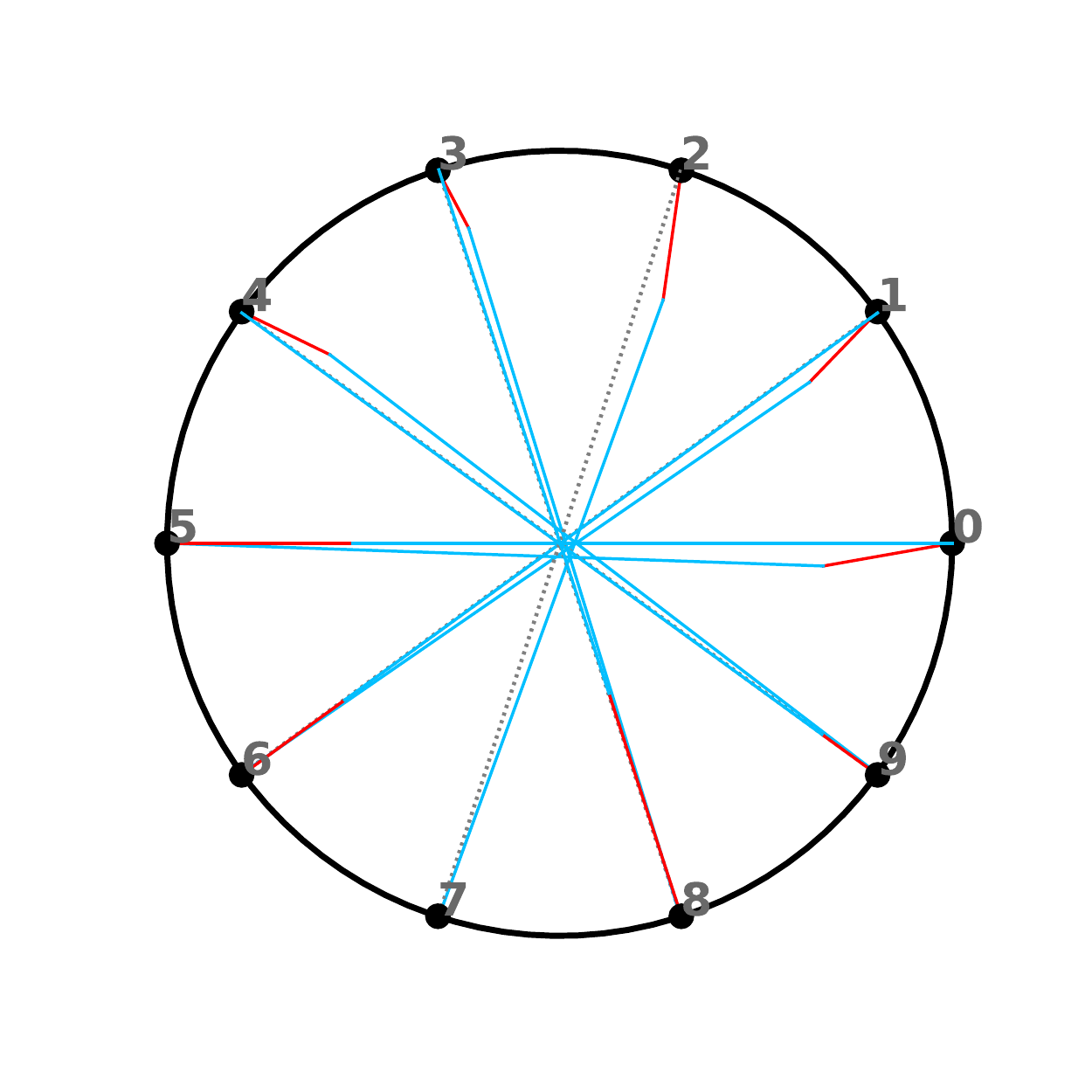}\label{cp2}}\  
	\subfloat[RCP-10-2]{\includegraphics[width=0.2\textwidth]{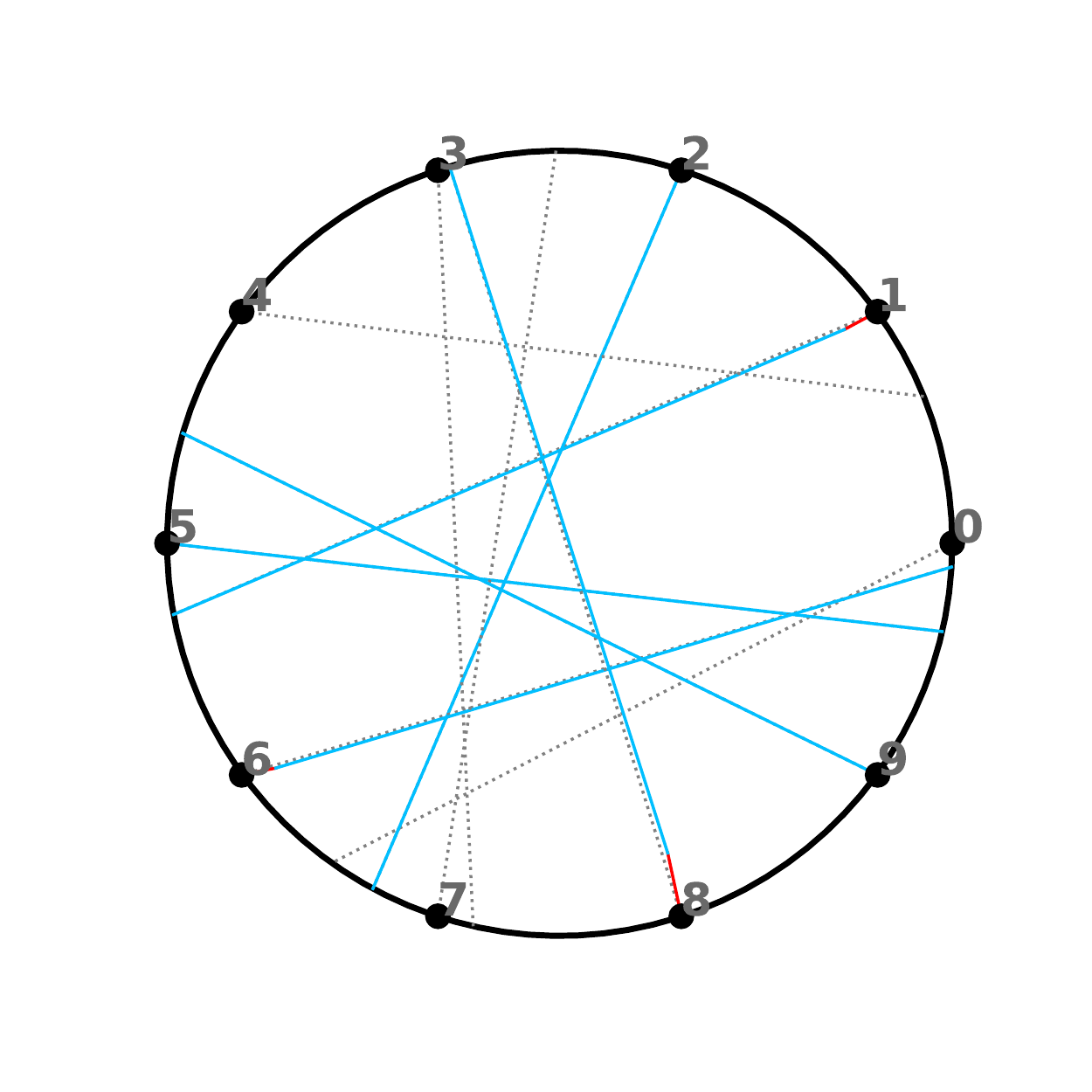}\label{rcp3}}\
	\subfloat[RCP-20-2]{\includegraphics[width=0.2\textwidth]{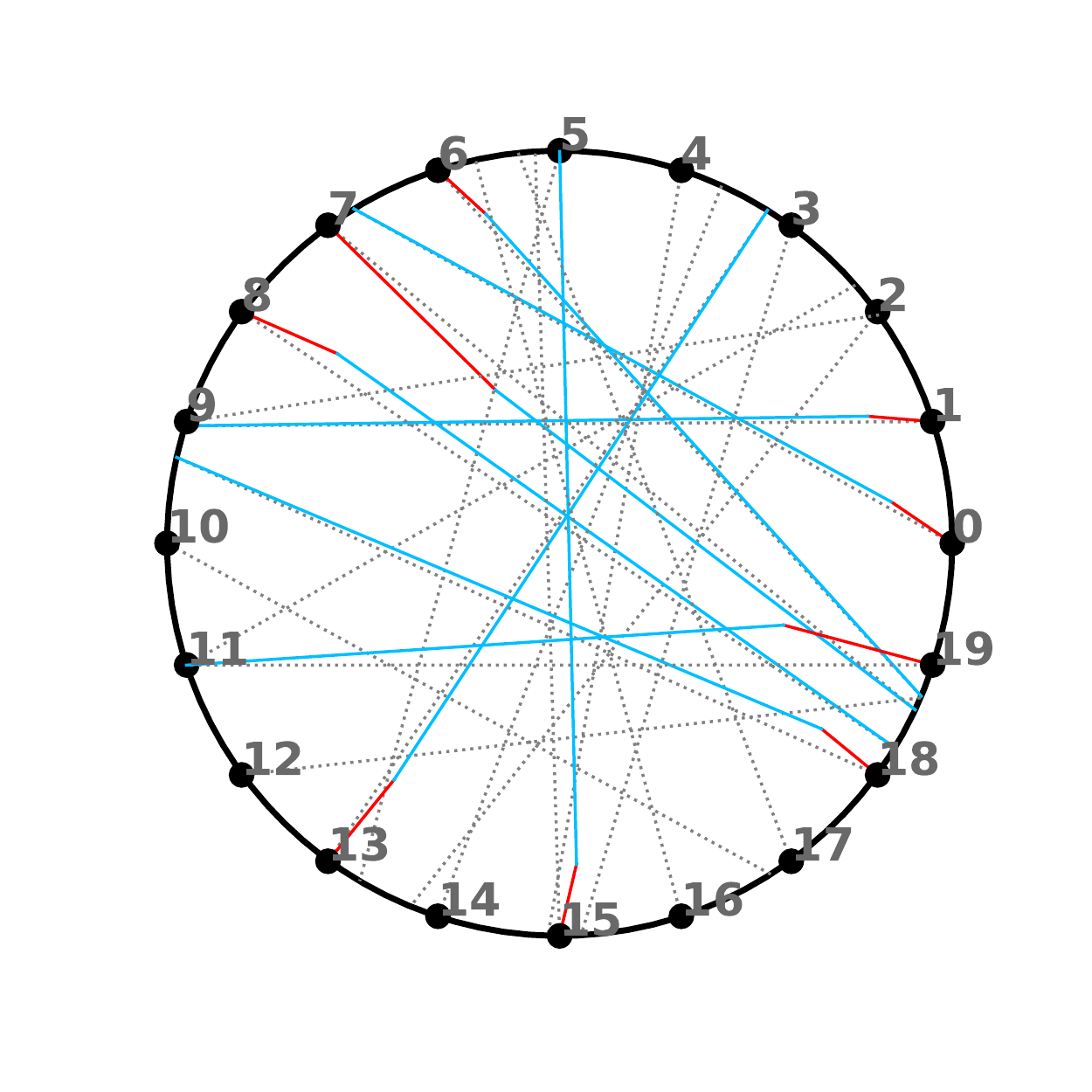}\label{cp3}}\  
	\subfloat[RCP-30-2]{\includegraphics[width=0.2\textwidth]{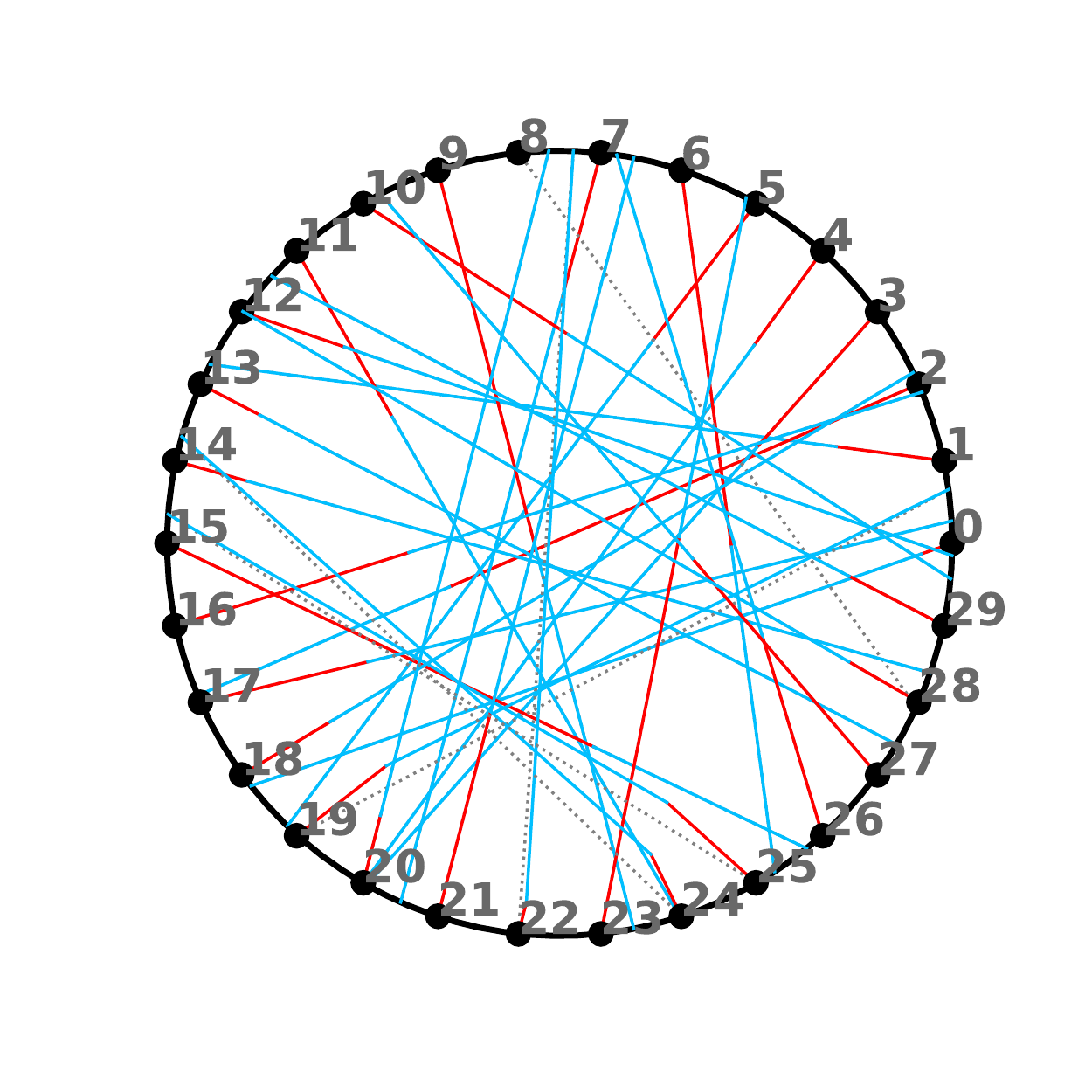}\label{rcp5}}\
	\subfloat[CP-15]{\includegraphics[width=0.2\textwidth]{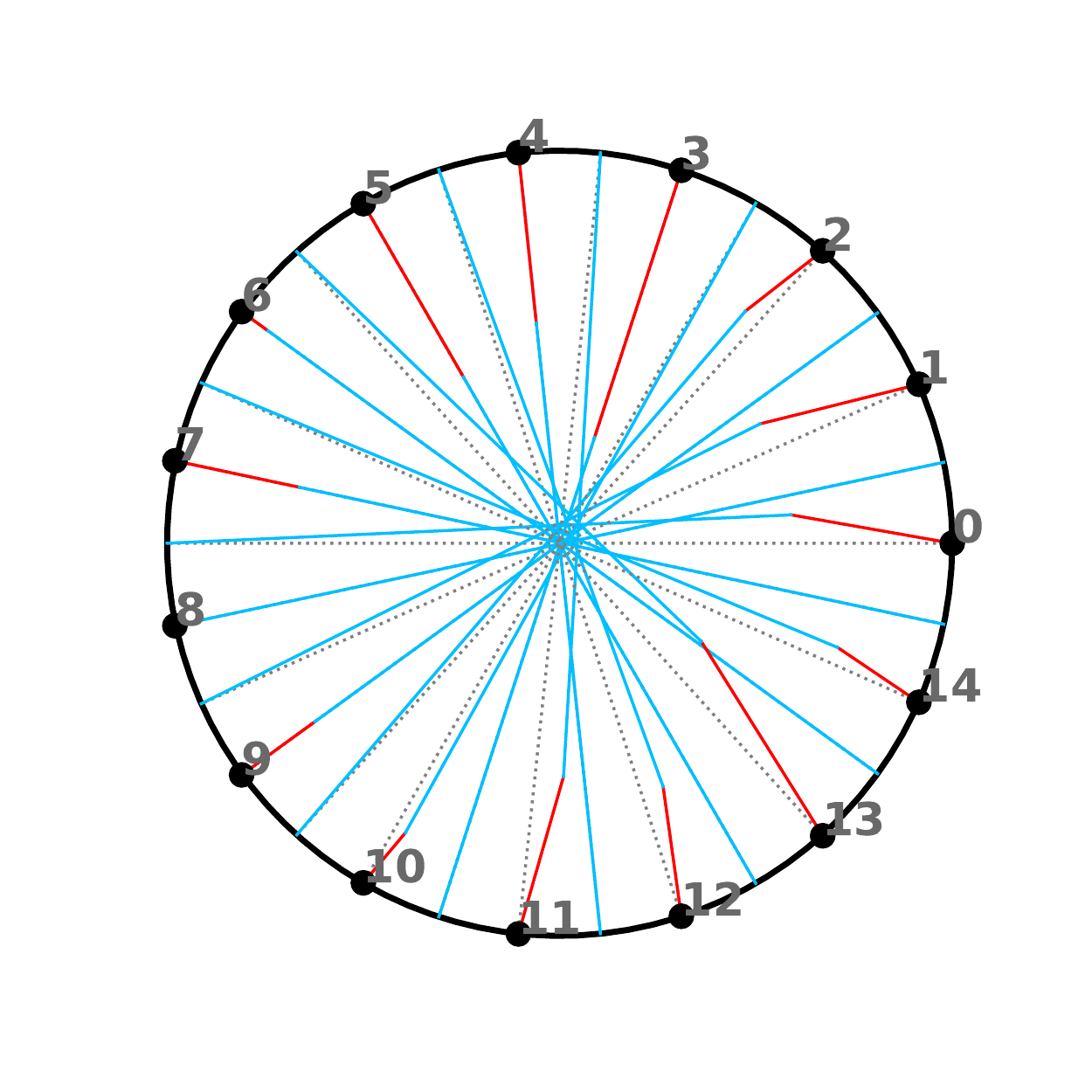}\label{cp4}}\
	\subfloat[RCP-10-3]{\includegraphics[width=0.2\textwidth]{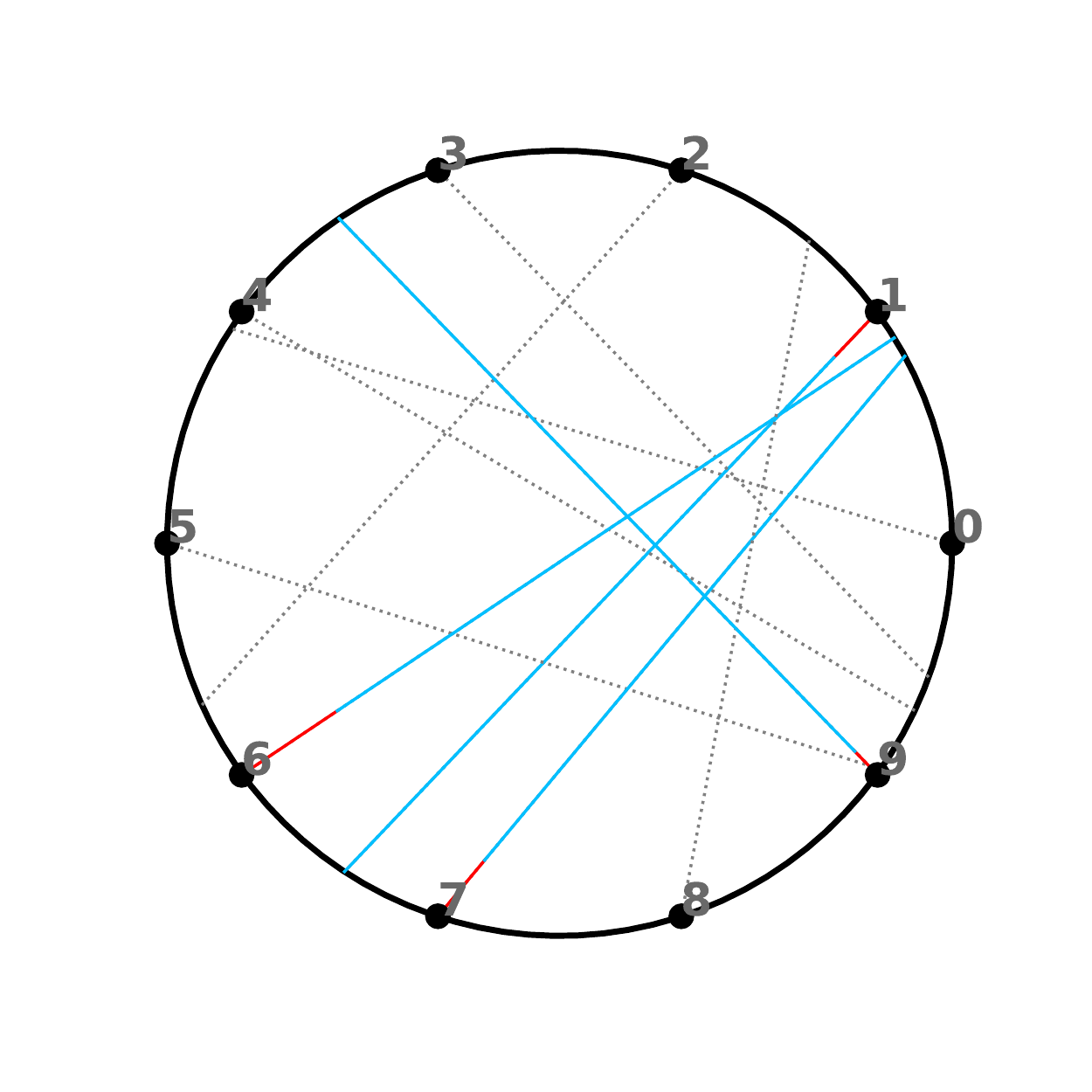}\label{rcp42}}\ 
	\subfloat[RCP-20-3]{\includegraphics[width=0.2\textwidth]{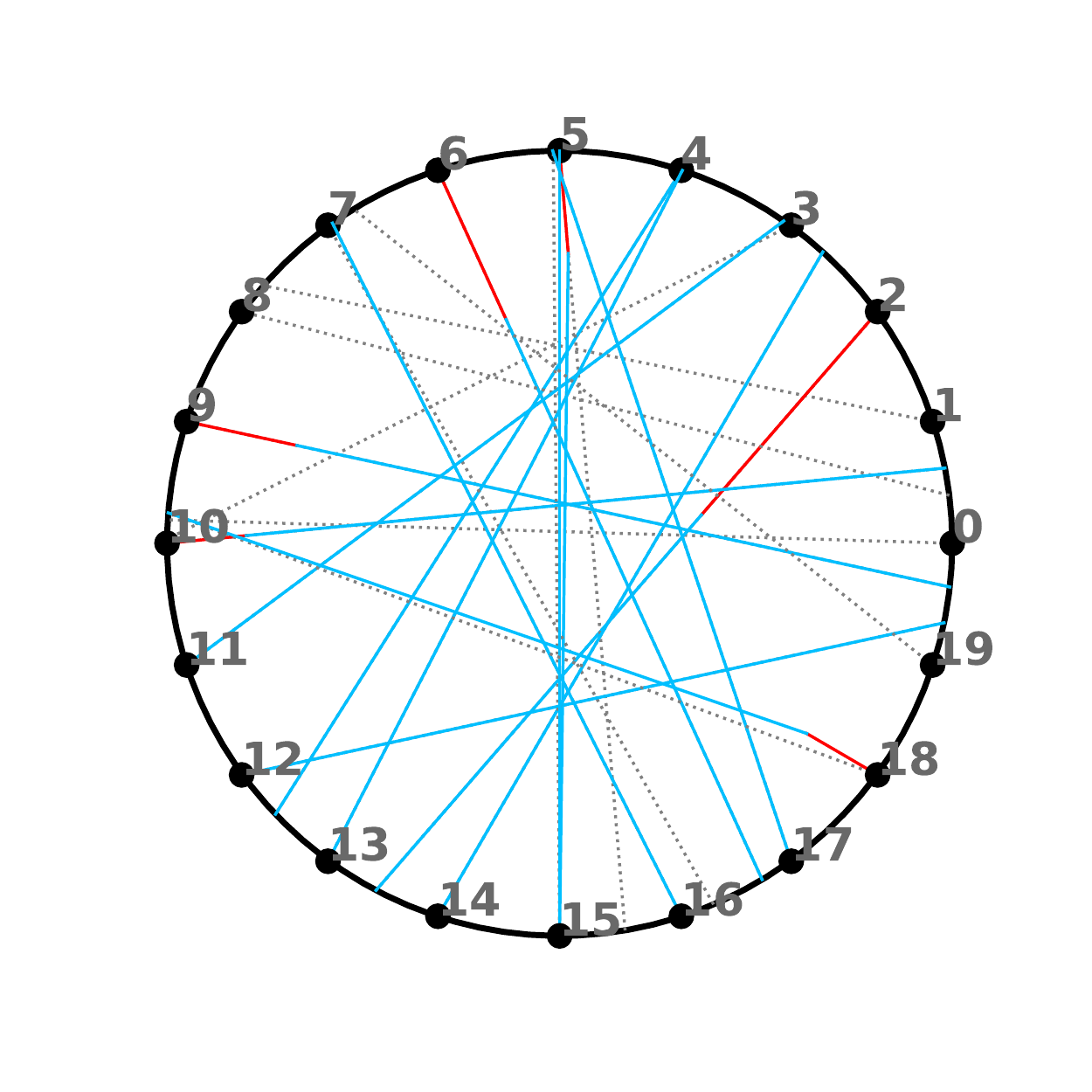}\label{rcp52}}\ 
	\subfloat[RCP-30-3]{\includegraphics[width=0.2\textwidth]{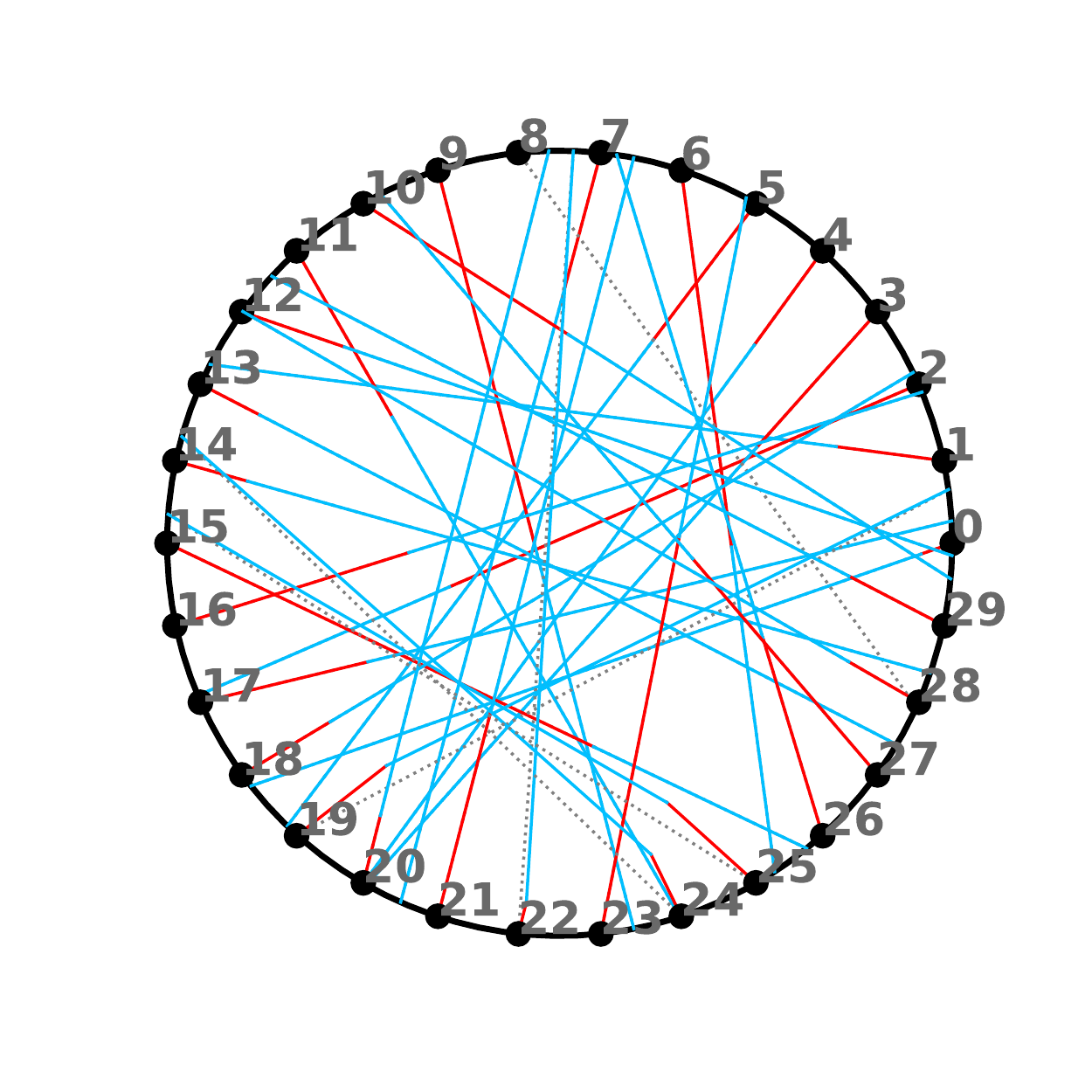}\label{rcp6}}
	\caption{Illustration of the two-stage algorithm with Exact-Recovery. Dashed grey lines represent aircraft initial trajectories. Red lines represent aircraft collision avoidance trajectories obtained via Model \ref{mod:dm}. Blue lines represent recovery trajectories obtained via Model \ref{mod:onestepdetailed}.}
\end{figure*}\label{ilus} 

\subsection{Performance of the two-stage algorithm}

We examine the performance of the two variants of the two-stage algorithm, one where the trajectory recovery stage is solved using Model \ref{mod:onestepdetailed} (Exact Recovery -- ER), and the other using Algorithm \ref{algo:greedy} (Greedy Recovery -- GR). Table \ref{tab:CP} summarizes the results for CP instances. In the header $|\mathcal{A}|$ is the number of aircraft, $n_c$ as the number of conflicts, \textit{Obj.} is the objective function and \textit{Time(s)} is the runtime in seconds. We also report the minimum recovery time among all aircraft $\min\limits_{i \in \mathcal{A}} t_i$, the average recovery time $\frac{1}{|A|}\sum\limits_{i \in \mathcal{A}}{t_i}$, and the maximum recovery time $\max\limits_{i \in \mathcal{A}} t_i$. \textit{Gap\%} is the relative gap difference between the objective values of ER and GR.

In terms of performance, we find that the runtime of Models \ref{mod:dm} and \ref{mod:onestepdetailed} increase exponentially with the number of aircraft, highlighted the challenging nature of the problems. Comparing trajectory recovery approaches, GR scales very efficiently in terms of runtime compared to ER is able to find some optimal solutions. Overall the average recovery time of GR is comparable to that of ER. In some cases, ER is able to override the decision at the first stage by recovering aircraft initially deviated at $t=0$, which means that such aircraft do not need to perform any avoidance maneuver. 

\begin{table*}[!h]
	\resizebox{2\columnwidth}{!}{%
		\begin{tabular}{ll ll lllll lllll l }
			\toprule
			&&\multicolumn{2}{l}{Avoidance} & \multicolumn{5}{l}{Exact-Recovery} & \multicolumn{5}{l}{Greedy-Recovery}\\
			\cmidrule(l){3-4} \cmidrule(l){5-9} \cmidrule(l){10-14}
			$|A|$ & $n_c$ & Obj. & Time (s) & Obj. & Time (s) & $\min\limits_{i \in \mathcal{A}} t_i$  & $\frac{1}{|A|}\sum\limits_{i \in \mathcal{A}}{t_i}$ & $\max\limits_{i \in \mathcal{A}} t_i$  & Obj. & Time (s) & $\min\limits_{i \in \mathcal{A}} t_i$ & $\frac{1}{|A|}\sum\limits_{i \in \mathcal{A}}{t_i}$ & $\max\limits_{i \in \mathcal{A}} t_i$ & Gap(\%) \\
			\midrule
4 & 6 & 4.89E-02 & 0.02 & 3.25E-03 & 0.24 & 0.00 & 1.25 & 2.00 & 3.25E-03 & 0.06 & 0.00 & 1.00 & 2.00 & 0.06 \\ 
5 & 10 & 3.08E-03 & 0.03 & 9.76E-05 & 0.35 & 0.00 & 2.00 & 4.00 & 3.19E-04 & 0.08 & 0.00 & 2.40 & 5.00 & 69.41 \\
6 & 15 & 7.34E-02 & 0.06 & 4.89E-03 & 0.47 & 0.00 & 2.17 & 4.00 & 8.14E-03 & 0.10 & 0.00 & 2.17 & 5.00 & 39.9\\
7 & 21 & 2.86E-02 & 0.52 & 2.01E-03 & 0.87 & 0.00 & 2.71 & 7.00 & 3.84E-03 & 0.12 & 1.00 & 3.86 & 7.00 & 47.5\\
8 & 28 & 9.82E-02 & 0.11 & 8.23E-03 & 0.97 & 2.00 & 3.50 & 6.00 & 8.23E-03 & 0.14 & 0.00 & 2.38 & 5.00 & 0.02\\
9 & 36 & 7.72E-02 & 26.9 & 7.81E-03 & 1.52 & 1.00 & 3.33 & 6.00 & 7.81E-03 & 0.22 & 0.00 & 1.78 & 5.00 & 0.04\\
10 & 45 & 1.23E-01 & 5.48 & 1.32E-02 & 2.03 & 0.00 & 3.50 & 7.00 & 1.88E-02 & 0.25 & 0.00 & 3.30 & 10 & 29.8\\
11 & 55 & 1.05E-01 & 50.4 & 1.44E-02 & 17.1 & 2.00 & 4.73 & 8.00 & 2.17E-02 & 0.29 & 0.00 & 5.00 & 11.0 & 33.2\\
12 & 66 & 1.49E-01 & 39.6 & 1.91E-02 & 9.74 & 0.00 & 4.25 & 9.00 & 2.62E-02 & 0.26 & 0.00 & 3.67 & 10.0 & 27.0\\
13 & 78 & 1.52E-01 & 6.93 & 2.53E-02 & 301 & 0.00 & 4.62 & 8.00 & 3.43E-02 & 0.38 & 0.00 & 5.31 & 9.00 & 26.2\\
14 & 91 & 1.74E-01 & 130 & 2.57E-02 & 40.0 & 0.00 & 4.36 & 8.00 & 2.56E-02 & 0.36 & 0.00 & 3.21 & 8.00 & 0.01\\
15 & 105 & 1.98E-01 & 200 & 2.81E-02 & 225 & 0.00 & 4.07 & 8.00 & 3.70E-02 & 0.49 & 0.00 & 4.27 & 7.00 & 23.8\\
			\bottomrule
	\end{tabular}}
	\caption{Results on the Circle Problem with 4 to 15 aircraft.}
	\label{tab:CP}
\end{table*}

\begin{figure}[!h]
	\centering
	\includegraphics[width=1\linewidth]{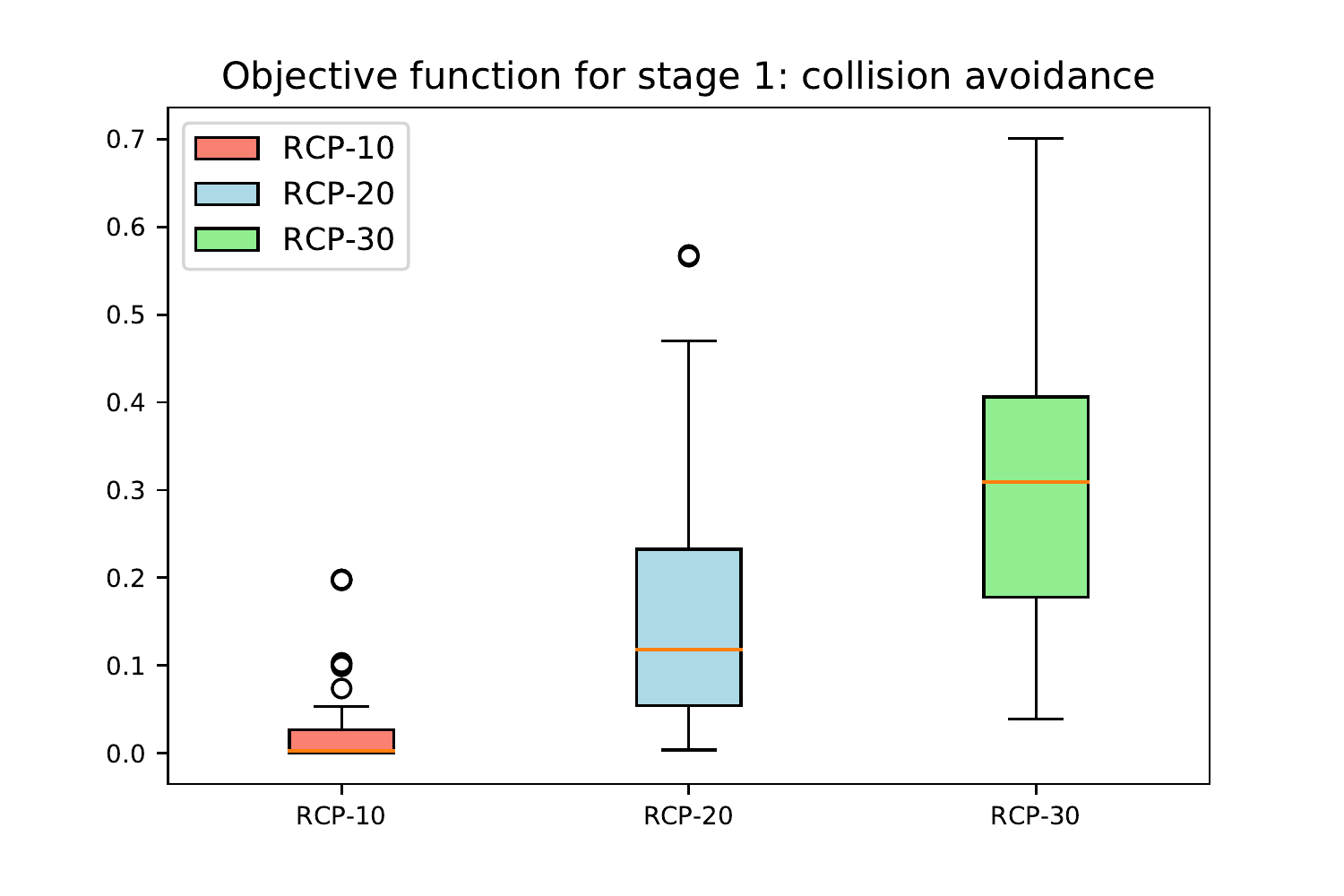}
	\caption{Stage 1: collision avoidance. Distribution of objective function values for Model \ref{mod:dm}. Each boxplot represent 100 instances of each instance size (10, 20 and 30 aircraft).}
	\label{fig:obj}
\end{figure}

\begin{figure}[!h]
	\centering
	\includegraphics[width=1\linewidth]{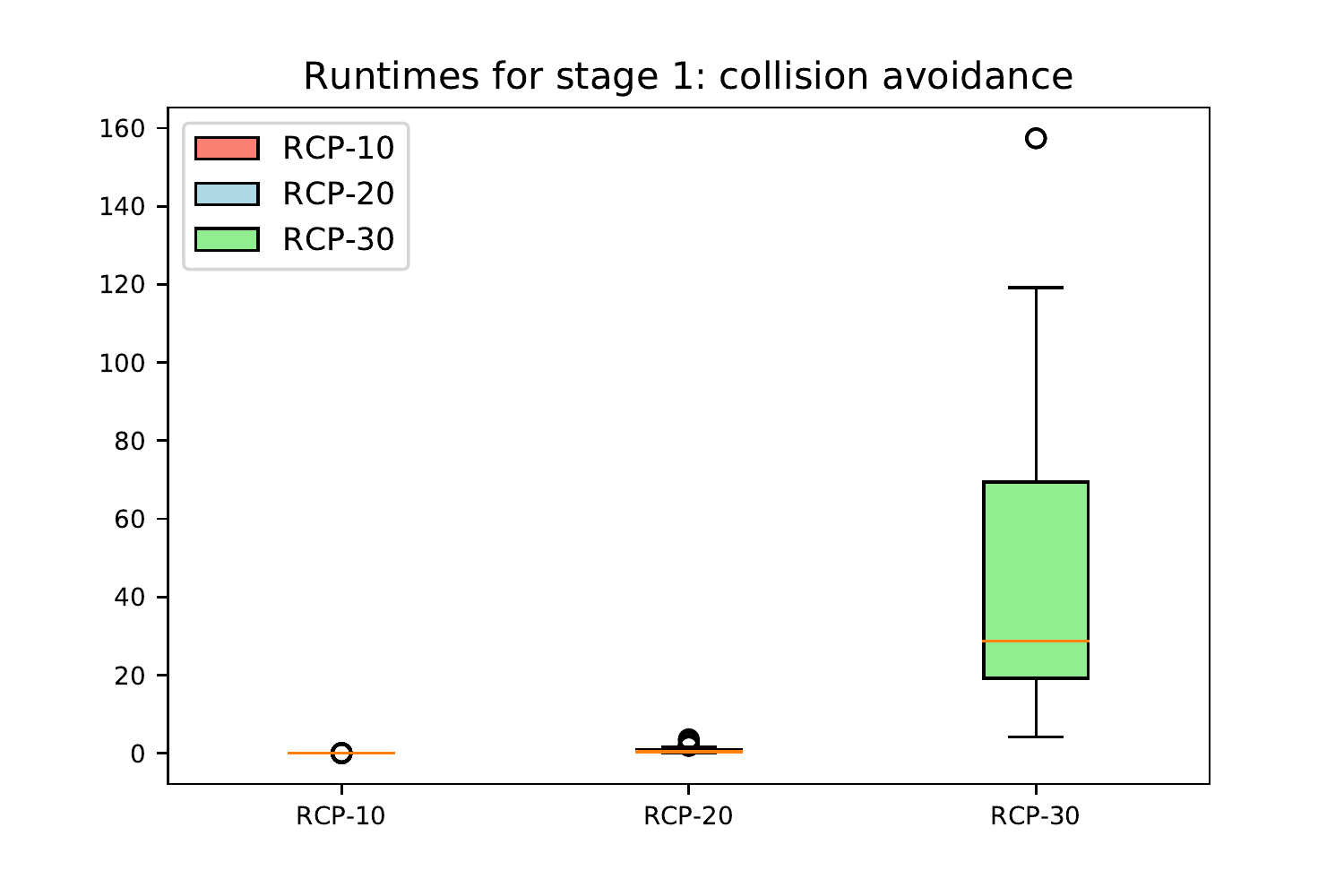}
	\caption{Stage 1: Collision avoidance. Distribution of runtimes for Model \ref{mod:dm}. Each boxplot represent 100 instances of each instance size (10, 20 and 30 aircraft).}
	\label{fig:runtime}
\end{figure}

\begin{figure}[!h]
	\centering
	\includegraphics[width=1\linewidth]{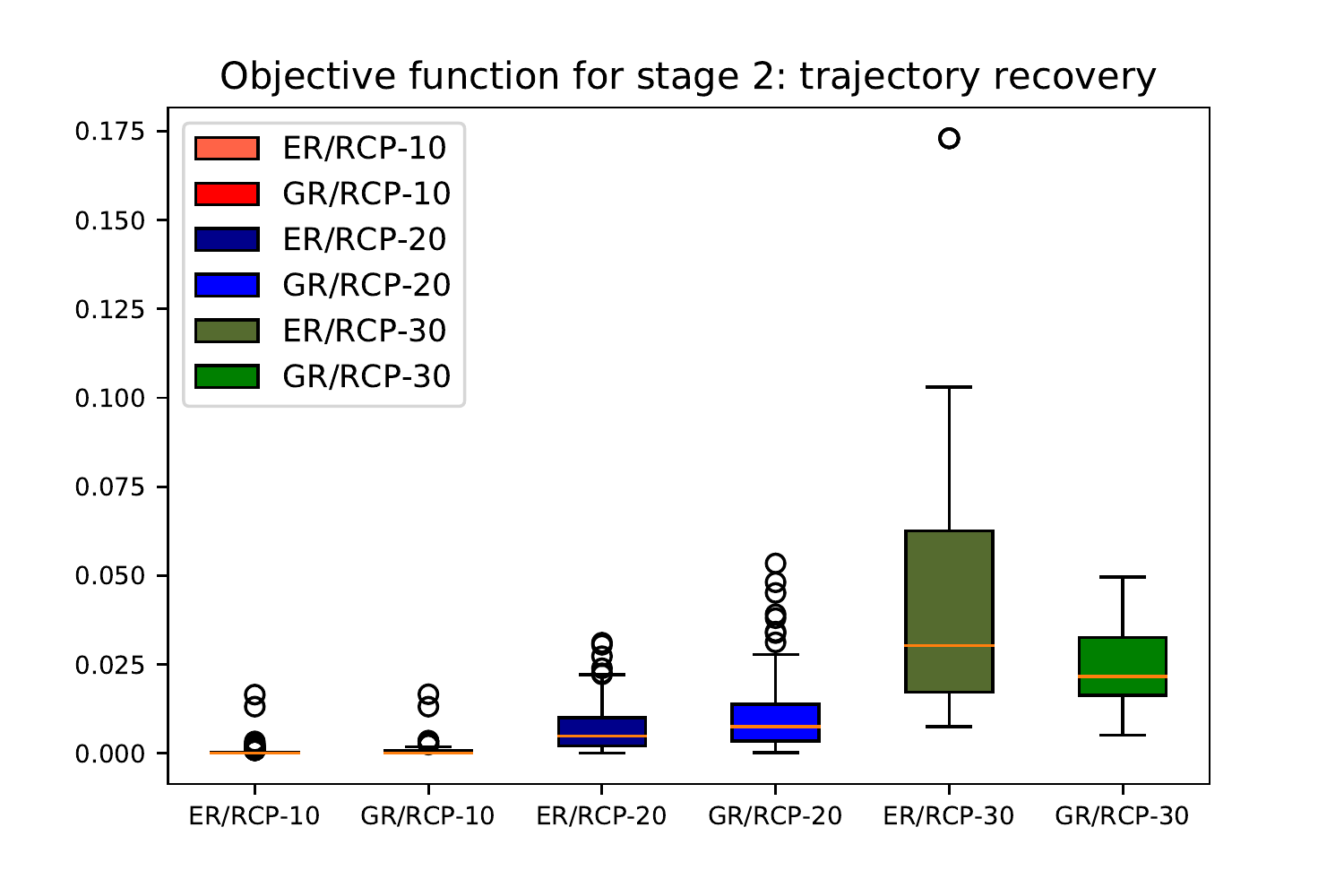}
	\caption{Stage 2: Trajectory recovery. Distribution of objective function values for ER (Model \ref{mod:onestepdetailed}) and GR (Algorithm \ref{algo:greedy}). Each boxplot represent 100 instances of each instance size (10, 20 and 30 aircraft).}
	\label{fig:obj2}
\end{figure}

\begin{figure}[!h]
	\centering
	\includegraphics[width=1\linewidth]{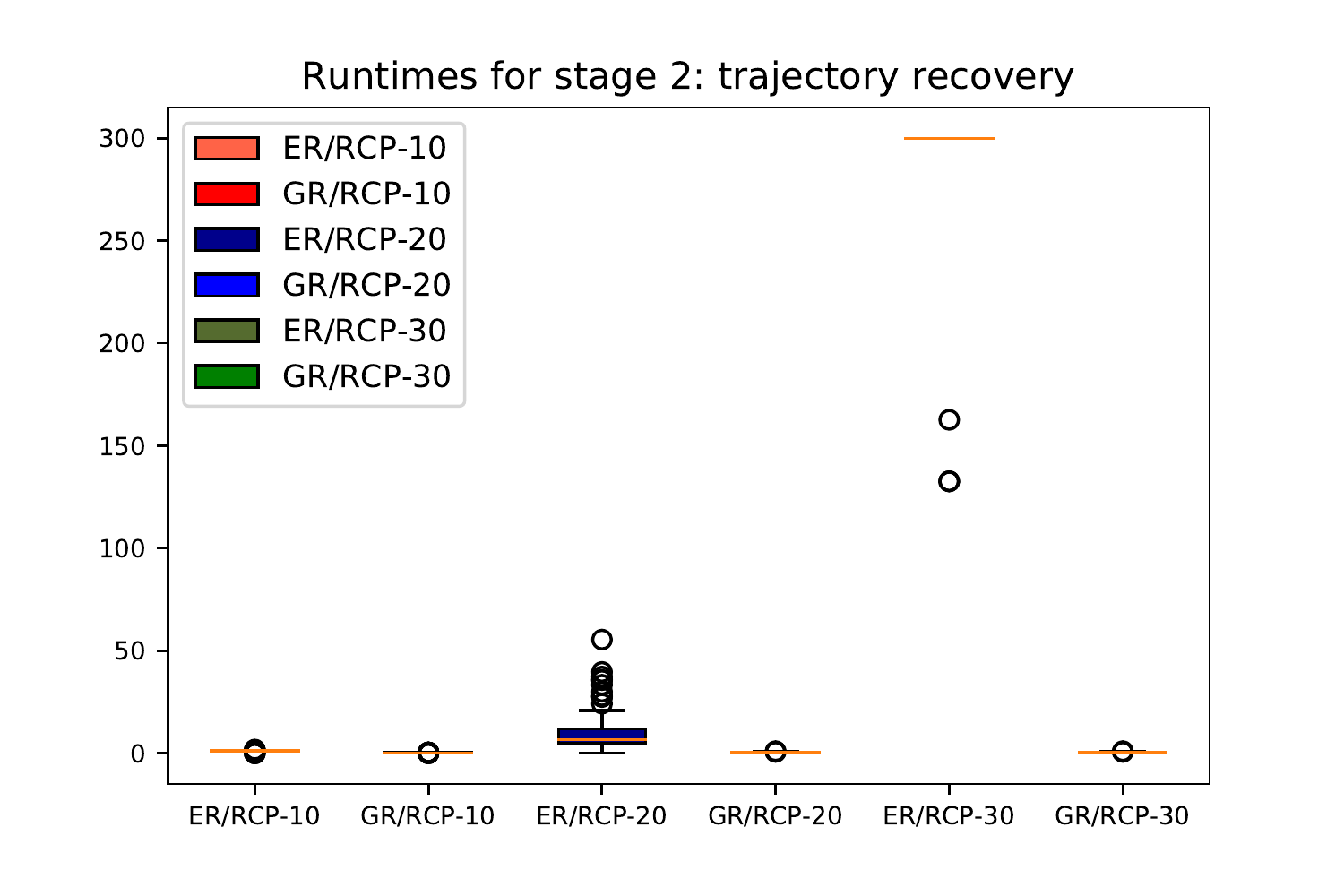}
	\caption{Stage 2: Trajectory recovery. Distribution of runtimes for ER (Model \ref{mod:onestepdetailed}) and GR (Algorithm \ref{algo:greedy}). Each boxplot represent 100 instances of each instance size (10, 20 and 30 aircraft).}
	\label{fig:runtime2}
\end{figure}

We report the performance of the two-stage algorithms for RCP instances in Figures \ref{fig:obj}-\ref{fig:runtime2}. For stage 1, we observe that the objective function and its variance increase super-linearly with the number of aircraft (see Figure \ref{fig:obj}). All first stage problems are solved within less than 3 minutes, with 10- and 20-aircraft instances requiring less than 10 seconds (see Figure \ref{fig:runtime}).

For stage 2, the performance of ER and GR on RCP instances is found to be comparable in terms of objective function values (see Figure \ref{fig:obj2}). For 10- and 20-aircraft problems, GR is marginally sub-optimal. For 30-aircraft problems, we find that when using a time limit of 5 minutes for ER most instances time out (see Figure \ref{fig:runtime2}), and the feasible solution returned by ER is often less competitive than the one provided by GR. This highlights the potential of the proposed greedy heuristic for real-time decision support. 

\section{Conclusion and perspectives}

We proposed a new two-stage algorithm for aircraft conflict resolution with trajectory recovery. Our approach decomposes the problem by first solving initial conflicts using heading and speed control; before identifying the optimal time for recovery in a second stage. The first stage is solved via a MIQP which extends an existing formulation to discretized heading control. Two novel solution methods have been proposed for the trajectory recovery state: an exact MILP formulation and a greedy heuristic algorithm. Both approaches use discrete time, thus the trajectory recovery problem consists of finding optimal aircraft conflict-free recovery times among a set of possible alternatives. Our objective is to minimize deviation in the avoidance stage and recovery time by accounting for the avoidance cost in the recovery stage. The performance of the proposed formulations is tested on benchmark problems for conflict resolution. We find that the proposed two-stage algorithm is able to solve instances with up to 30 aircraft  using the exact collision avoidance and trajectory recovery models. The proposed greedy algorithm for trajectory recovery is found to systematically construct feasible solutions of near optimal quality on random test problems thus offering a competitive alternative to scale-up the approach. 

Further testing is required to fully assess the impact of the control variables onto solution quality. Notably, quantifying the impact of first stage formulations with non-discretized heading changes and second stage formulation in continuous time is critical to evaluate the cost of maneuver discretization. Future research will also be focused on improving the interplay between both stages of the algorithm. Emphasis will be placed on developing stochastic formulations capable of accounting for the expected cost of trajectory recovery at the collision avoidance stage. Further, future work will investigate the design of robust formulations to account for the uncertainty of aircraft trajectory prediction when planning. 

\bibliographystyle{elsarticle-harv}
\bibliography{biblio}

\end{document}